\documentclass[a4paper]{article}
\usepackage[left=1in,top=1in,right=1in,bottom=1in,nohead]{geometry}
\usepackage{authblk}
\usepackage{blindtext}
\usepackage[T1]{fontenc}
\usepackage[utf8]{inputenc}
\usepackage{lmodern}
\usepackage{lipsum}
\usepackage[english]{babel}
\usepackage{csquotes}
\usepackage[english]{babel}
\usepackage{amsfonts}
\usepackage{amssymb}
\usepackage{amsmath}
\usepackage{amsthm}
\usepackage{eucal}
\usepackage{mathrsfs}
\allowdisplaybreaks[0]
\usepackage{bm}
\usepackage{graphicx}
\usepackage{caption}
\usepackage{subcaption}
\usepackage{algorithm,algorithmic}
\usepackage{framed}
\usepackage{diagbox}
\usepackage[titletoc,title]{appendix}
\usepackage{tikz}
\usepackage{color}
\usepackage{hyperref}
\usepackage{cleveref}
\usepackage{autonum}
\title{Overcomplete}
\date{}

\usepackage{mathrsfs}

\newcommand{\R}{{\mathbb R}}

\newcommand{\mM}{{\mathsf M}}
\newcommand{\mC}{{\mathsf C}}
\newcommand{\mR}{{\mathsf R}}
\newcommand{\mT}{{\mathsf T}}
\newcommand{\mL}{{\mathsf L}}
\newcommand{\mS}{{\mathsf S}}
\newcommand{\mA}{{\mathsf A}}

\newcommand{\mD}{{\mathsf D}}

\newcommand{\mI}{{\mathsf I}}
\newcommand{\mQ}{{\mathsf Q}}

\newcommand{\mO}{{\mathsf O}}
\newcommand{\mSigma}{{\mathsf \Sigma}}

\DeclareGraphicsExtensions{.pdf,.png,.jpg}

\newtheorem{theorem}{Theorem}

\usepackage{enumitem}
\setlist[enumerate]{leftmargin=.5in}
\setlist[itemize]{leftmargin=.5in}

\title{Sparsity promoting hybrid solvers for hierarchical Bayesian inverse problems}

\begin{document}

\title{Sparsity promoting hybrid solvers for hierarchical Bayesian inverse problems}

\author{Daniela Calvetti$^{\rm a}$, Monica Pragliola$^{\rm b}$, Erkki Somersalo$^{\rm a}$\\\vspace{6pt} 
$^{a}${\em{Department of Mathematics, Applied Mathematics and Statistics, Case  Western Reserve University}}\\
$^{b}${\em{Department of Mathematics, University of Bologna}}}
\maketitle

\begin{abstract}
	The recovery of sparse generative models from few noisy measurements is an important and challenging problem. Many deterministic algorithms rely on some form of $\ell_1$-$\ell_2$ minimization to combine the computational convenience of the $\ell_2$ penalty and the sparsity promotion of the $\ell_1$. It was recently shown within the Bayesian framework that sparsity promotion and computational efficiency can be attained with hierarchical models with conditionally Gaussian priors and gamma hyperpriors. The related Gibbs energy function is a convex functional and its minimizer, which is the MAP estimate of the posterior, can be computed efficiently with the globally convergent Iterated Alternating Sequential (IAS) algorithm \cite{CSS}. Generalization of the hyperpriors for these sparsity promoting hierarchical models to generalized gamma family yield either globally convex Gibbs energy functionals, or can exhibit local convexity for some choices for the hyperparameters. \cite{CPrSS}. The main problem in computing the MAP solution for greedy hyperpriors that strongly promote sparsity is the presence of local minima.  To overcome the premature stopping at a spurious local minimizer, we propose two hybrid algorithms that first exploit the global convergence associated with gamma hyperpriors to arrive in a neighborhood of the unique minimizer, then adopt a generalized gamma hyperprior that promote sparsity more strongly.  The performance of the two algorithms is illustrated with computed examples.
\end{abstract}

\section{Introduction}
The recovery of a sparse vector from noisy indirect observations continues to be an active research topic. After the groundbreaking work on compressed sensing and its connections to sparsity-promoting regularization methods \cite{BDE,CandesRombergTao,CandesTao,Donoho2006,DET} and the $\ell_1$-penalty in particular, the interest in sparse recovery has been revived by dictionary learning methods in data science, where the goal is to match an observed vector with few dictionary entries in a huge data base \cite{Kreutz,Mairal,Tariyal}.
The connections between regularization methods and penalty functionals on one hand, and Bayesian inference techniques on the other, have been thoroughly investigated  \cite{CSbook,KSbook,CSWIRES}, and families of priors that promote sparsity have been identified in the Bayesian framework. 

Sparsity is a qualitative rather than quantitative trait because in general the size of the support and its location cannot be specified in advance. While there is a wealth of different priors that promote sparsity, the results may differ significantly depending on the cost for non-vanishing entries. In the classical regularization setting, this is well illustrated by the different $\ell_p$-penalties, with $p\leq 1$. Penalty functionals with $0\leq p <1$ tend to promote sparsity more strongly than $p=1$.  The convexity of the objective function for the latter, and the results on the exact recovery of sparse generative models under suitable conditions have contributed to the popularity of $\ell_1$ regularization for sparse problems, while the presence of local minima has damped the enthusiasm for penalties with $p<1$.  

Analogous consideration hold in the Bayesian framework for sparsity promoting hierarchical prior models, with generalized gamma hyperpriors. Recently \cite{CSS,CPrSS}, it has been shown that the Maximum A Posteriori (MAP) iterative sparse reconstruction algorithm is particularly well suited for heavily underdetermined but large scale problems (see, e.g., \cite{CPPSV} for an application). The Iterative Alternating Sequential  algorithm (IAS) is based on hierarchical Bayesian models, and uses sparsity promoting hyperpriors selected from a family of generalized gamma distributions. As pointed out in \cite{CPrSS}, some choices of the hyperparameters yield algorithms that are closely related, e.g., to the $\ell_p$-penalization methods. Moreover, the convexity properties of the objective function also depend on the parameter choice, as does the convergence rate of algorithms for computing the MAP estimate. Our aim is to combine the properties of generalized gamma hyperpriors to design robust and computationally efficient methods for sparse recovery from few noisy observations. More specifically, we propose hybrid algorithms for the MAP computation: a  gamma hyperprior guides the approximate solution towards the unique minimizer of the objective functions at the beginning, and subsequently a greedier hyperprior is employed to promote sparsity more  strongly.  We focus on two such hybrid algorithms, that we refer to as local and global because of the different strategy to switch hyperpriors. In the local version, the hyperprior is changed componentwise, guaranteeing local convexity, while in the global version, the hyperprior is changed for all components. In addition to analyzing the convergence properties of each approach, we provide a criterion for ensuring that the a priori beliefs are consistent with the two different hyperpriors. The performance of the algorithms is assessed in the light of computed examples.   

\section{Hierarchical Bayesian models}
\label{sec:bayes}

In this section we introduce the Iterative Alternating Sequential (IAS) algorithm for the MAP computation, and review some of its key properties: additional details can be found in \cite{CPPSV,CSS,CPrSS}.

Consider the linear observation model with additive Gaussian noise,
\begin{equation}
b = \mA x + e\, , \quad e \sim \mathcal{N}(0,\mSigma),
\label{eq:lin_model}
\end{equation}
where $\mA\in \R^{m \times n}$, with $m<n$, is a known ill-conditioned matrix describing the forward model, $x\in\R^n$ is the unknown of interest, $\mSigma \in \R^{n\times n}$ is the symmetric positive definite covariance matrix of the noise. We remark that by letting $\mA' = \mS\mA$ and $b'=\mS b$, where $\mS$ is the Cholesky factor of the precision matrix,
$\mSigma^{-1} = \mS^{\mT}\mS$, we can assume the noise to be white, i.e. $\mSigma=\mI$, hence, the likelihood probability density function (pdf) of $b$ with given $x$ takes the form
\begin{equation}
\pi_{b\mid x}(b\mid x) \propto \exp\bigg(-\frac{1}{2}\| \mA x - b \|^{2}\bigg).
\label{eq:like}
\end{equation}
We are interested in estimating $x$ from the observed measurements in $b$ under the a priori assumption that $x$ is sparse, that is, $\|x\|_0 ={\rm card}({\rm supp}(x))\ll n$. In general, the approach can be generalized to cases where the unknown of interest itself is not sparse, but admits a sparse representation in some dictionary,  by making the coefficients of the representation the unknown of primary interest. 

To encode the sparsity belief in the prior model, we begin by considering a componentwise Gaussian prior model, 
\begin{equation}
x_j \sim \mathcal{N}(0,\theta_j)\, ,\quad  \theta_j > 0\, , \quad 1 \leq j \leq n,
\label{prior1}
\end{equation}
or equivalently, 
\begin{equation}
x \sim \mathcal{N}(0,\mD_{\theta})\, ,\mD_{\theta} = \mathrm{diag}(\theta_1,\ldots,\theta_n) \in \R^{n\times n},
\label{prior2}
\end{equation}
where the variances of the individual components are not known. The conditional prior density of $x$ given $\theta$ is of the form
\begin{equation}
\pi_{x\mid\theta} (x\mid \theta) \propto \frac{1}{\prod_{j=1}^{n}\sqrt{\theta_j}}\exp\bigg({-\frac{1}{2}\| \mD_{\theta}^{-1/2} x \|^{2}}\bigg)
=\exp\bigg({-\frac{1}{2}\| \mD_{\theta}^{-1/2} x \|^{2}} -\frac 12\sum_{j=1}^n \log\theta_j\bigg),
\label{prior3}
\end{equation}
and following the Bayesian paradigm that all unknowns are modeled as random variables, the a priori belief about $\theta$ is encoded in a \emph{hyperprior} pdf $\pi_{\Theta}(\theta)$.  The price to pay for this hierarchical prior model is that we need to estimate not only $x$ but also $\theta$ based on data in terms of the joint posterior distribution of $(x,\theta)$ conditioned on $b$, 
\begin{equation}
\pi_{x,\theta\mid b}(x,\theta\mid b) \propto  \underbrace{\pi_{x\mid \theta}(x\mid\theta)\,\pi_\theta(\theta)}_{\pi_{x,\theta}(x,\theta)} \,\pi_{b\mid x}(b\mid x).
\label{post}
\end{equation}
A way to promote sparse solutions is to choose a hyperprior $\pi_\theta$  that favors small values of $\theta$ but allows occasional  large outliers.
A family with these properties, thoroughly investigated in \cite{CPrSS}, is that of the generalized gamma distributions,
\begin{equation}
\pi_{\theta}(\theta) = 
\pi_{\theta}(\theta \mid r,\beta,\vartheta) = \frac{|r|^n}{\Gamma(\beta)^n} \prod_{j=1}^{n}\frac 1{\vartheta_j} \left(\frac{\theta_j}{\vartheta_{j}}\right)^{r\beta - 1} \exp\bigg(-\left(\frac{{\theta_j}}{{\vartheta_{j}}}\right)^r\;\bigg)\,,
\label{eq:hyper}
\end{equation}
where $r\in \R\setminus\{0\}$, $\beta >0$, $\vartheta_j>0$.

The  MAP estimate, of the posterior pdf model (\ref{post}) is also the minimizer of the negative logarithm of the posterior pdf,
\begin{equation}
(x^{*},\theta^{*}) =\arg\min_{x,\theta}\, \big\{ -\log \pi_{x,\theta\mid b}(x,\theta\mid b) \,{=:}\,\mathcal{F}(x,\theta)\big\}.
\label{eq:map}
\end{equation}
The objective function ${\mathcal  F}(x,\theta)$ can be written as 
\begin{eqnarray}
\nonumber
&&\mathcal{F}(x,\theta) = \mathcal{F}(x,\theta \mid r,\vartheta,\beta)\\
&=&\quad \lefteqn{ \overbrace{\phantom{\frac 12 \|b - \mA  x\| ^2 + \frac 12 \sum_{j=1}^n\frac{x_j^2}{\theta_j} }}^{(a)}}
{\displaystyle\frac 12} \|b - \mA  x\|^2+
\underbrace{ \frac 12\sum_{j=1}^n\frac{x_j^2}{\theta_j}
	-\bigg(r\beta - \frac 32\bigg)
	\sum_{j=1}^n\log\frac{\theta_j}{\vartheta_j} +\sum_{j=1}^n\left(\frac{\theta_j}{\vartheta_j}\right)^r}_{(b)}, \label{split}
\end{eqnarray}
to emphasize that  only the terms in (a) depend on $x$, and only those in (b) depend on $\theta$. These observations play a key role for the design of a computationally efficient algorithm for computing the MAP estimate.  We start by recalling the IAS algorithm for the solution of problem (\ref{eq:map});  see \cite{CPSV,CPPSV,CSS} for further details and for a comprehensive study of the effect of the choice of hyperparameters $(r,\beta,\vartheta)$ on the promotion of sparsity and the properties of the objective function.   

\subsection{IAS algorithm}
Given the initial value $\theta^0$, each step of the IAS for problem (\ref{eq:map}) consists of the two updates,
\[
\theta^t \rightarrow x^{t+1}\rightarrow \theta^{t+1},\quad t\geq 0,
\] 
where 
\begin{equation}
\label{z upd pb}
x^{t+1} = \arg\min_{x}\,\left\{\mathcal{F}(x,\theta^t)\right\}\,,\quad \theta^{t+1} = \arg\min_{\theta}\,\left\{\mathcal{F}(x^{t+1},\theta)\right\}\,.
\end{equation}
Due to the particular form of the objective function (\ref{split}), each step comprises first the computation of the minimizer of (a) with respect to $x$ keeping $\theta$ fixed, then the minimizer of (b) with respect to $\theta$ with the updated value of $x$ fixed. While this procedure is remarkably similar to the Alternating Direction Method of Multipliers (ADMM) \cite{Boyd},  there are some fundamental differences. In fact, while in ADMM,  the alternating structure is achieved via an artificial partial decoupling of the fidelity term and the penalty term by introducing auxiliary variables, in IAS the partial decoupling is automatic, with the common term of (a) and (b) being the link between the two minimization tasks. Moreover, both minimization tasks are relatively simple with an exact condition for the minimizer.  For some choices of hyperparameters, IAS has been shown to be globally at least linearly convergent \cite{CPSV,CSS}. In the following, we review some of the computational details of the IAS algorithm that are particularly relevant for the proposed hybrid schemes.

\paragraph{Update of $x$}
The update of $x$ given $\theta$  by minimizing part (a) in (\ref{split}) reduces to solution of a quadratic minimization problem, i.e., 
\[
x^{t+1} = \arg\min_{x}\,\left\{
\|\mA x - b \|^2 + \|\mD_{\theta}^{-1/2}x \|^2
\right\}\,,\;\theta=\theta^{t}\,,
\]
thus $x^{t+1}$ is the least squares solution of the linear system
\begin{equation}
\label{exact LSQ}
\left[\begin{array}{c}
\mA\\
\mD_{\theta}^{-1/2}
\end{array}\right] x = 
\left[\begin{array}{c}
b\\
0
\end{array}\right]\,.
\end{equation}
The solution of (\ref{exact LSQ}) can be approximated by solving a reduced problem via the Conjugate Gradient for Least Squares (CGLS) algorithm \cite{CPSV2}, often without any real loss of information in the solution \cite{CPrSS}. Introduce the change of variables,
\[
\mD_{\theta}^{-1/2}x = w, 
\]  
which corresponds to a whitening of the conditional prior,  and reformulate the linear system (\ref{exact LSQ}) in terms of $w$ as
\begin{equation}
\label{exact LSQ 2}
\left[\begin{array}{c}
\mA_{\theta}\\
\mI
\end{array}\right] w = 
\left[\begin{array}{c}
b\\
0
\end{array}\right]\,,\;\mA_{\theta} = \mA\mD_{\theta}^{1/2}.
\end{equation}
If the matrix $\mA$ is ill-conditioned, the linear system
\begin{equation}
\label{x upd}
\mA_{\theta}w = b\,,\quad x= \mD_{\theta}^{1/2} w,
\end{equation}
can be solved approximately through a CGLS iteration with an early stopping  criterion \cite{CPrSS}. More precisely, denote the $k$th Krylov subspace corresponding to the above system by
\[
{\mathscr K}_k =  {\mathscr K}_k(\mA_\theta^\mT b,\mA^\mT_\theta\mA) = {\rm span}\left\{\left(\mA_\theta^\mT\mA_\theta\right)^j \mA_\theta b \mid 0\leq j\leq k-1\right\}.
\]  
Define the Reduced Krylov Subspace (RKS) solution as 
\[ 
w_k = {\rm argmin}\left\{ \|b - \mA_\theta w\| \mid w \in {\mathscr K}_k\right\},
\]
where $k$ is the first index satisfying
\[
\|b - \mA_\theta w_{k+1}\|\leq \sqrt{m}, \quad G(w_{k+1}) >\tau G(w_k),
\]
where $\tau>1$, $\varepsilon = \tau-1>0$ is a small safeguard parameter, and the functional $G$, given by 
\[
G(w) = \|b - \mA_\theta w\|^2 + \|w\|^2,
\]
is the objective function approximately minimized by the surrogate reduced model. 

\paragraph{Update of $\theta$} It follows from the independence of the components that the first order optimality condition that needs to be satisfied by the updated $\theta$ can be imposed componentwise. Setting the partial derivatives of (b) in ({\ref{split}) with respect to $\theta_j$ equal to zero, we find that $\theta_j$ must satisfy
	\begin{equation}
	\label{get xi}
	-\frac{1}{2} \frac{x_j^2}{\theta_j^2}-\left(r \beta -\frac{3}{2}\right)\frac1{\theta_j} + r \frac{\theta_j^{r-1}}{\vartheta_j^r} = 0\,,\quad x=x^{t+1}\,.
	\end{equation}
	While for some values of $r$, notably $r=\pm 1$, (\ref{get xi}) admits an analytic solution, in general we need to solve it numerically. It was shown in \cite{CPrSS} that after the changes of variables  $\theta_j = \vartheta_j \xi_j$, $x_j = \sqrt{\vartheta_j} z_j$, we may write $\xi_j = \varphi(|z_j|)$, and  via implicit differentiation, the function $\varphi$ satisfies the initial value problem
	\begin{equation}\label{get xi 2}
	\varphi'(z) = \frac{2z \varphi(z)}{2r^2 \varphi(z)^{r+1}+z^2} , \quad \varphi(0) = \left(\frac{\eta}{r}\right)^{1/r},
	\end{equation}
	therefore the updated $\theta_j$ can be computed by a numerical time integrator.
	
	We conclude this section with the main results on selecting the model parameters $(r,\beta,\vartheta)$.  The values of the parameters $r$ and $\beta$ affect how strongly the sparsity of the solution is promoted and determine the convexity of the objective function, while the value of $\vartheta_j$ can be related to the sensitivity of the data to $x_j$. Recall that for a linear model $b = \mA x + \varepsilon$, a classical measure of the sensitivity of the data $b$ to the component $x_j$ is 
	$ \|\mA e_j\| $, where $e_j\in\R^n$ is the canonical $j$th Cartesian unit vector. It was proven recently \cite{CPPSV, CSS, CPrSS} that under rather natural conditions, a judicious choice of the parameter $\vartheta$ is 
	\[
	\vartheta_j = \frac{C}{\|\mA e_j\|^2},
	\]
	where the constant $C>0$ is related to the expected sparsity of of the solution and to an estimate of the signal-to-noise ratio (SNR). Due to the connection with sensitivity, we this choice of $\vartheta$ is referred to as sensitivity scaling.

	\section{Hybrid IAS algorithms}
	\label{sec:hyb hyp}
	In this section, we will propose a hybrid version of the IAS algorithm in which the hypermodel in the generalized gamma family is updated componentwise as the iteration proceeds. The following theorem, see \cite{CPrSS} for details, summarizes how the values of the hyperparameters $r$ and $\beta$ affect the convexity properties of the functional $\mathcal{F}$ .
	
	\begin{theorem}\label{th:convexity}
		Let $\beta>0$ and $r\neq 0$, and let $\mathcal{F}(x,\theta)$ be the objective function for the minimization problem in (\ref{eq:map}).
		\begin{itemize}
			\item[(a)] 
			If $r\geq 1$ and $\eta = r\beta  -3/2>0$, the function $\mathcal{F}(x,\theta)$ is globally convex.
			\item[(b)] 
			If $0<r<1$ and $\eta = r\beta  -3/2>0$, or, if $r<0$ and $\beta>0$,
			the function $\mathcal{F}(x,\theta)$ is convex provided that
			\begin{equation} \label{th bar}
			\theta_j < \overline{\theta}=\vartheta_j \left(\frac{\eta}{r|r-1|}\right)^{1/r}.
			\end{equation}
		\end{itemize}
	\end{theorem}
	
	As far as the computation of the MAP estimate is concerned, the global convexity of the objective function when $r\geq 1$ is very convenient, although there are several reasons for considering other choices of $r$ that yield hierarchical priors that promote sparsity more strongly. It has been observed that, by and large, the further the objective function is from being globally convex, the stronger the sparsity of the minimizer is promoted. We review below some recent results, see  \cite{CHPS,CSS}, relating generalized gamma hyperpriors and classical sparsity-promoting priors.
	
	Let
	\begin{eqnarray*}
		{\mathcal P}(x,\theta\mid r,\beta,\vartheta) &=& \frac 12\sum_{j=1}^n\frac{x_j^2}{\theta_j}
		-\bigg(r\beta - \frac 32\bigg)
		\sum_{j=1}^n\log\frac{\theta_j}{\vartheta_j} -\sum_{j=1}^n\left(\frac{\theta_j}{\vartheta_j}\right)^r \\
		&=& \sum_{j=1}^n p(x_j,\theta_j\mid r,\beta,\vartheta_j) 
	\end{eqnarray*}
	denote the penalty term (b) in the objective function, and express the IAS updating formula (\ref{get xi 2}) for $\theta_j$ as a function of $x_j$
	\begin{equation}\label{def g}
	g_j(x_j) =  \theta_j = \vartheta_j \varphi\left(\frac{|x_j|}{\sqrt{\vartheta_j}}\right).
	\end{equation}
	It has been shown in \cite{CSS} that for $r=1$, as $\eta\to 0+$ the penalty function ${\mathcal P}(x,\theta,1,3/2+\eta,\vartheta)$ approaches a weighted $\ell_1$-penalty in the sense that 
	\[
	\lim_{\eta\to 0+}  {\mathcal P}(x, g(x) \mid 1,\frac{3}{2}+\eta,\vartheta) = \sqrt{2} \sum_{j=1}^n \frac{|x_j|}{\sqrt{\vartheta_j}}
	\]
	and, moreover, the corresponding minimizer $x^*$  found by the IAS algorithm converges to scaled $\ell_1$ regularized solution. 
	
	More generally as shown in \cite{CPrSS}, by choosing $r\beta = 3/2$, the penalty function coincides with the weighted $\ell_p$-norm, with $p =2r/(r+1)$,
	\[
	{\mathcal P}\big(x, g(x)\mid r,\frac{3}{2r},\vartheta\big)  = C_r\sum_{j=1}^n \frac{|x_j|^p}{\sqrt{\vartheta_j}^p}, \quad C_r = \frac{r+1}{(2r)^{r/(r+1)}}.
	\] 
	While this result holds in general, for $0<r<1$ and $\beta = 3/2r$, the model corresponds to $\ell_p$ penalties with $0<p<1$, which are known to promote  strongly the sparsity of the solution.
	
	For the inverse gamma hypermodel, corresponding to $r=-1$, the penalty term approaches the Student distribution, a prominently fat tailed distribution favoring large outliers, and leading to a greedy algorithm that strongly promotes sparsity \cite{CPrSS}. The main problem with the  lack of global convexity of the objective function is that optimization-based algorithms for the MAP computation may stop at a spurious local minimizer.
	
	In this work, we propose two modifications to the IAS algorithm that take advantage of the global convexity of the objective function corresponding to the gamma hyperprior ($r=1$), and of the stronger sparsity promotion of hierarchical models with $r<1$ whose associated objective functions are only locally convex. In both proposed algorithms, the gamma hyperprior is used initially to drive the IAS iterates towards the unique minimizer of the globally convex objective function, then switched to greedier hypermodel. 
	The two different algorithms are referred to as {\em local} and {\em global} hybrid models. In the local hybrid, the hyperprior is changed componentwise as soon as the corresponding variance falls inside the convexity region of the second model, while in the global model, the hyperprior is changed for all components after a given number of IAS steps.  Next we present the details relative to the two hybrid schemes.

	\subsection{Local hybrid IAS}
	We write the objective function ${\mathcal F}(x,\theta\mid r,\vartheta,\beta)$ with the given model parameters $(r,\beta,\vartheta)$  as
	\[
	{\mathcal F}(x,\theta\mid r,\vartheta,\beta) = \|b - \mA x\|^2 + \sum_{j=1}^n p(x_j,\theta_j \mid r,\vartheta_j,\beta),
	\]
	where
	\[
	p(x_j,\theta_j \mid r,\vartheta_j,\beta) = \frac 12 \frac{x_j^2}{\theta_j}  -\bigg(r \beta  - \frac{3}{2}\bigg) \log \frac{\theta_j}{\vartheta_j} + \left(\frac{\theta_j}{\vartheta_j}\right)^r.
	\]  
	Unlike in the standard IAS algorithm, where the parameters $r$, $\beta$ and $\vartheta$ are kept fixed, the local hybrid algorithm updates the parameters for those component pairs $(x_j,\theta_j)$ that satisfies the convexity criterion in Theorem~\ref{th:convexity} for the second hypermodel. 
	
	More precisely, consider two hypermodels  with parameters $(r^{(1)},\beta^{(1)},\vartheta^{(1)})$  and  $(r^{(2)},\beta^{(2)},\vartheta^{(2)})$, with $r^{(2)}<1\leq r^{(1)}$, $r^{(2)}\neq 0$, referred to as ${\mathscr M}_1$ and ${\mathscr M}_2$, respectively, and start the IAS algorithm with the model ${\mathscr M}_1$.

	Let $(x,\theta) = (x^t,\theta^t)$ denote the IAS iterate after  $t$ steps. For each component  $x_j$ of $x$, we compute the $\theta_j$ update corresponding to model ${\mathscr M}_2$, 
	\[
	\theta_j^{(2)} = g(x_j\mid r^{(2)},\beta^{(j)},\vartheta^{(2)}_j) = g^{(2)}(x_j).
	\]
	If 
	\begin{equation}\label{switch cond}
	\theta_j^{(2)} < \overline \theta_j = \vartheta^{(2)}_j \left(\frac{\eta^{(2)}}{r^{(2)}|r^{(2)}-1|}\right)^{1/r^{(2)}},
	\end{equation} 
	we update $\theta_j$ switching to ${\mathscr M}_2$, otherwise we continue with  ${\mathscr M}_1$. Observe that since the function $g^{(2)}$ is strictly increasing for $x_j>0$, we may write the above condition in terms of $x_j$,
	\[
	|x_j| < \left[g^{(2)}\right]^{-1}(\overline\theta_j) = \overline x_j.
	\] 
	Let $I\subset \{1,2,\ldots,n\}$ denote an index set such that
	\[
	j\in I \mbox{ if and only if $|x_j|<\overline x_j$,}
	\] 
	and by $I^c$ its complement. Define the local hybrid objective function,
	\begin{eqnarray*}
		{\mathcal F}(x,\theta\mid I) &=& \|b - \mA x\|^2 + \sum_{j\in I^c} p(x_j,\theta_j\mid r^{(1)},\vartheta_j^{(1)},\beta^{(1)}) +  \\
		&&\phantom{XXXXXXXXXX}
		\sum_{j \in I} p (x_j,\theta_j\mid r^{(2)},\vartheta_j^{(2)},\beta^{(2)}).
	\end{eqnarray*}
	whose convexity can be guaranteed by a bound constraint
	\begin{equation}
	\label{bound hyb}
	|x_j|< \overline x_j \mbox{ for $j\in I$}.
	\end{equation}
	It was shown in \cite{CPrSS} that to add a bound constraint to the IAS algorithm it suffices to project the updated vector $x$ onto the feasible set.
	The selection of the hyperparameter $\vartheta^{(j)}$, $j=1,2$, deserves some attention. For ${\mathscr M}_1$, the value of $\vartheta^{(1)}$ can be decided by taking sensitivity analysis into consideration, as suggested in  \cite{CPrSS}. We assign the value of $\vartheta^{(2)}$ based on the following consideration: {\em If $x_j=0$, the corresponding variance $\theta_j$ should be the same regardless of the choice of the hypermodel, and should reflect the expected variance of a background signal.} We recall that if $x_j=0$, the updating of $\theta_j$ in the IAS algorithm according to (\ref{get xi}) yields
	\[
	g(0\mid r,\beta,\vartheta_j) =  \vartheta_j \left(\frac{\eta}r\right)^{1/r},\quad \eta = r\beta - 3/2, 
	\] 
	and in order for the two models to agree, it suffices to set 
	\[
	\vartheta^{(2)}_j =  \left(\frac{\eta^{(1)}}{r^{(1)}}\right)^{1/r^{(1)}}\left(\frac{r^{(2)}}{\eta^{(2)}}\right)^{1/r^{(2)}}\vartheta_j^{(1)}.
	\]
	We are now ready to summarize the proposed local hybrid IAS scheme in algorithmic form. Here we assume that $x\in\R^n$ itself is sparse; suitable adjustments need to be made when the sparsity assumption concerns the increments.

	\begin{algorithm}
		\caption{Local Hybrid IAS}
		\label{alg:1}
		\vspace{0.2cm}
		{\renewcommand{\arraystretch}{1.2}
			\renewcommand{\tabcolsep}{0.0cm}
			\vspace{-0.08cm}
			\begin{tabular}{ll}
				\textbf{inputs}:      & Noisy data $\,b\in \R^m$,
				\vspace{0.04cm} \\
				&linear forward operator $\mA\in\R^{m\times n}$, noise covariance matrix $\mSigma\in\R^{m\times m}$	\vspace{0.04cm} \\
				&hyperparameters $(r^{(1)},\beta^{(1)},\vartheta^{(1)})$,  $(r^{(2)},\beta^{(2)},\vartheta^{(2)})$
				\vspace{0.04cm} \\
				\textbf{output}:$\;\;$     & estimated signal and variance $\,x^*,\theta^* \in \R^n$ \vspace{0.2cm} \\
			\end{tabular}
		}
		\vspace{0.1cm}
		{\renewcommand{\arraystretch}{1.2}
			\renewcommand{\tabcolsep}{0.0cm}
			\begin{tabular}{rcll}
				1. & $\quad$ & \multicolumn{2}{l}{\textbf{initialize:}$\;\;$ 
					set $t=0$, $\theta^t = \vartheta^{(1)}$, $I=\emptyset$ } \vspace{0.05cm}\\
				2. && \multicolumn{2}{l}{\textbf{for} $\;$ \textit{t = 0, 1, 2, $\, \ldots \,$ until convergence $\:$} \textbf{do}:} \vspace{0.1cm}\\
				3. && $\quad\;\;$ update $x^{t+1}$ by solving (\ref{x upd}) & 
				\vspace{0.05cm} \\
				4. && $\quad\;\;$ project components $x^{t+1}_j$, $j\in I$, to $[-\overline x,\overline x]$ & 
				\vspace{0.05cm} \\	
				5. && \multicolumn{2}{l}{$\quad\;\;$ \textbf{for}$\;$ \textit{j = 1,$\, \ldots \,,n\:$}   
				} \vspace{0.05cm} \\
				6. && $\quad\quad\quad\;\;$ \textbf{if}$\;$ $\theta_j \geq \overline{\theta}$  & 
				\vspace{0.05cm} \\
				7. && $\quad\quad\quad\quad\;\;$ update $\theta_j^{t+1} = g(x_j^{t+1}\mid r^{(1)},\beta^{(1)},\vartheta_j^{(1)})$ & 
				\vspace{0.05cm} \\
				8. &&$\quad\quad\quad\;\;$ \textbf{else}$\;$ & 
				\vspace{0.05cm} \\
				9. && $\quad\quad\quad\quad\;\;$ update $\theta_j^{t+1} = g(x_j^{t+1}\mid r^{(2)},\beta^{(2)},\vartheta_j^{(2)})$& 
				\vspace{0.05cm} \\
				10. && $\quad\quad\quad\quad\;\;$ update $I = I\cup\{j\}$& 
				\vspace{0.05cm} \\
				11.  &&$\quad\quad\quad\;\;$ \textbf{endif}$\;$ & 
				\vspace{0.05cm} \\
				
				12. && \multicolumn{2}{l}{\textbf{end$\;$for}} \vspace{0.09cm} \\
				13. && \multicolumn{2}{l}{$x^* = x^{t+1},\,\theta^{*}=\theta^{t+1}$}
			\end{tabular}
		}
	\end{algorithm}
	
	Before discussing a modification of the above algorithm, a comment on the projection on convexity interval  (step 4) is of order. The projection step is included in the algorithm to ensure that the index set $I$ of components being updated using the hypermodel ${\mathscr M}_2$ is monotonically increasing, which, in general, may not be automatically guaranteed. However, the numerical experiments show that the projection step in practice may not be necessary, and the bound constraint  $|x_j|< \overline x$ is not active.
	
	In \cite{CPrSS}, the stability of the convexity condition was briefly discussed in terms of the scaled (dimensionless) variables, $z_j = x_j/\sqrt{\vartheta_j^{(2)}}$, $\xi_j = \theta_j/\vartheta_j^{(2)}$.
	It was shown (see  Lemma 4.2 in \cite{CPrSS}) that if $\xi_j^t<\overline\xi$, then
	\[
	| z_j^{t+1}|\leq M\xi_j^t =M\varphi(|z_j^{t}|),
	\]
	where $\varphi$ is the IAS updating function of the scaled variable $\xi_j$ given the current $z_j$, and $M$ depends on the matrix $\mA$ and the data $b$. For example, in the case $r=1/2$, the convexity bound is
	$\overline \xi =  16\eta^2$.
	A natural question is whether, if $|z_j^t|<\varepsilon<\overline z$, where $\overline z=\varphi^{-1}(\overline\xi)$ is the convexity bound for the scaled variable $z_j$, it can be guaranteed that $|z_j^{t+1}|<\overline z$, or, equivalently,
	\begin{equation}\label{bound}
	\varphi(|z_j^{t+1}|)< \overline\xi = 16 \eta^2,
	\end{equation}
	where $ \eta = r^{(2)}\beta^{(2)} - 3/2$. In \cite{CPrSS} it was shown that
	\[
	\varphi(t) =  4 \eta^2 + \frac 1\eta t^2 + {\mathcal O}(t^4),
	\]
	therefore,
	\[
	| z_j^{t+1}| \leq M(4\eta^2 + \frac 1\eta \varepsilon^2).
	\] 
	To check if  (\ref{switch cond}) is satisfied up to fourth order terms, it suffices to have
	\[
	4\eta^2 + \frac 1\eta( M(4\eta^2 + \frac 1\eta \varepsilon^2))^2<16\eta^2,
	\] 
	or
	\[
	M^2(4\eta^2 + \frac  1\eta \varepsilon^2)^2< 12\eta^3,
	\]
	that is,
	\[
	\varepsilon^2< \frac{\sqrt{12}}M \eta^{3/2} - 4\eta^2.
	\]  
	The positivity of the right side can be guaranteed by choosing $\eta$ sufficiently small. While the above estimate is only approximate and qualitative, it conveys the idea that the stability of the convexity condition may depend on the forward model as well as on hyperparameter selection.

	\subsection{Global hybrid IAS}
	\label{sec:hyb ias}
	As the numerical experiments confirm, the gain from switching to the model ${\mathscr M}_2$ for components that are already in the convexity region of that model is not so much
	in enhancing, e.g., sudden discontinuities in the solution, but more on cleaning the background.
	An alternative approach is to relax the convexity requirement and use the global convexity of the first model to find a good starting point for the optimization with the second model, taking the risk of minimizing a non-convex objective function from an initial guess sufficiently close to the global minimum of the first model.  
	
	More specifically, we run first the IAS algorithm for a fixed number $\overline  t$ of iterations with model ${\mathscr M}_1$, whose conservative parameter choice guarantees convergence towards a global minimizer, and then switch to the less conservative hypermodel  ${\mathscr M}_2$, trading the global convexity for stronger sparsity promotion. We refer to this scheme as \emph{global} hybrid IAS, since the change of hyperprior is carried out at once for all the variances $\theta_j$.  unlike in the local version where only selected components followed the model ${\mathscr M}_2$. The computational details are summarized in Algorithm~\ref{alg:2}.

	\begin{algorithm}
		\caption{Global Hybrid IAS}
		\label{alg:2}
		\vspace{0.2cm}
		{\renewcommand{\arraystretch}{1.2}
			\renewcommand{\tabcolsep}{0.0cm}
			\vspace{-0.08cm}
			\begin{tabular}{ll}
				\textbf{inputs}:      & Noisy data $\,b\in \R^m$,
				\vspace{0.04cm} \\
				&linear forward operator $\mA\in\R^{m\times n}$, noise covariance matrix $\mSigma\in\R^{m\times m}$	\vspace{0.04cm} \\
				&hyperparameters $(r^{(1)},\beta^{(1)},\vartheta^{(1)})$,  $(r^{(2)},\beta^{(2)},\vartheta^{(2)})$ \\
				&integer $\overline t>0$ defining the switch point
				\vspace{0.04cm} \\
				\textbf{output}:$\;\;$     & estimated signal and variance $\,x^*,\theta^* \in \R^n$ \vspace{0.2cm} \\
			\end{tabular}
		}
		\vspace{0.1cm}
		{\renewcommand{\arraystretch}{1.2}
			\renewcommand{\tabcolsep}{0.0cm}
			\begin{tabular}{rcll}
				1. & $\quad$ & \multicolumn{2}{l}{\textbf{initialize:}$\;\;$ 
					set $\theta^{0} = \vartheta^{(1)}$ } \vspace{0.05cm}\\
				2. && \multicolumn{2}{l}{\textbf{for} $\;$ \textit{t = 0, 1, 2, $\, \ldots \,$ until convergence $\:$} \textbf{do}:} \vspace{0.1cm}\\
				3. && $\quad\;\;$ update $x^{t+1}$ by solving (\ref{x upd}) & 
				\vspace{0.05cm} \\	
				4. && \multicolumn{2}{l}{$\quad\;\;$ \textbf{for}$\;$ \textit{j = 1,$\, \ldots \,,n\:$}   
				} \vspace{0.05cm} \\
				5. && $\quad\quad\quad\;\;$ \textbf{if}$\;$ $t < \overline t$  & 
				\vspace{0.05cm} \\
				6. && $\quad\quad\quad\quad\;\;$ update $\theta_j^{t+1} = g(x_j^{t+1}\mid r^{(1)},\beta^{(1)},\vartheta_j^{(1)})$ & 
				\vspace{0.05cm} \\
				7. &&$\quad\quad\quad\;\;$ \textbf{else}$\;$ & 
				\vspace{0.05cm} \\
				8. && $\quad\quad\quad\quad\;\;$ update $\theta_j^{t+1} = g(x_j^{t+1}\mid r^{(2)},\beta^{(2)},\vartheta_j^{(2)})$& 
				\vspace{0.05cm} \\
				9.  &&$\quad\quad\quad\;\;$ \textbf{endif}$\;$ & 
				\vspace{0.05cm} \\
				
				10. && \multicolumn{2}{l}{\textbf{end$\;$for}} \vspace{0.09cm} \\
				11. && \multicolumn{2}{l}{$x^* = x^{(t+1)},\,\theta^{*}=\theta^{(t+1)}$}
			\end{tabular}
		}
	\end{algorithm}
	
	In the description of the Algorithm~\ref{alg:2}, the value $\overline t$ is given as input. Alternatively, one could run the model ${\mathscr M}_1$ until the variances $\theta$ stop changing significantly. Since in general we have little information of the nature of the minima of the objective function when $r<1$, a definitive automatic switching rule is not easy to justify.  We illustrate the performance of the algorithm on a few test cases in the section on computed examples.

	\section{IAS for sparse increments}
	\label{sec:ias incr}
	The IAS algorithm, and the hybrid versions of it,  assume that the unknown has a sparse representation, and $x$ is the vector of coefficients in this representation. In the case where the a priori sparsity belief is not about the signal $x$ but its increments, the IAS algorithms needs to be suitably adapted. In the  one dimensional case the changes are rather straightforward, while the treatment in the higher dimensional cases is more delicate. 
	
	In the one-dimensional case assume that the unknown is a piecewise constant signal in $\R$ characterized by few discontinuities. If $f(t)$ is the signal, $0\leq t\leq 1$ and  $x_j = f(jh)$, where $h = 1/n$ is the discretization parameter, we may express $x$ in terms of the increments,
	\[
	x_j = x_0 + \sum_{k=1}^j (x_k - x_{k-1}),  1\leq j\leq n,
	\]
	or, letting  $y_j = x_j-x_{j-1}$, as
	\[
	x = x_0 e_1 + \mL^{-1} y, 
	\]
	where
	\[
	\mL = \left[\begin{array}{rrrr} 1 & & & \\ -1 & 1 & & \\ &\ddots &\ddots & \\ & & -1 & 1\end{array}\right], \quad e_ 1 =\left[\begin{array}{c} 1 \\ 0 \\ \vdots \\ 0\end{array}\right].
	\] 
	Assuming, for simplicity that $x_0=0$, it follows from the invertibility of $\mL$ that we may reformulate the problem as estimating $y$ from the observation model
	\[
	b = \mA \mL^{-1} y +\varepsilon,
	\]
	with the a priori belief that $y$ is sparse.  We update $x$ in the IAS algorithm by computing
	\[
	y^{t+1} = {\rm argmin}\left\{ \frac 12 \|b - \mA \mL^{-1} y\|^2  + \frac 12\sum_{j=1}^n \frac{y_j^2}{\theta_j^{t}}\right\}, \quad x^{t+1} = \mL^{-1}y^{t+1},
	\]
	where the sparsity of $y$ is playing a role in the update of the $\theta$ 
	\[
	\theta^{t+1} = {\rm argmin}\left\{ 
	\frac 12\sum_{j=1}^n \frac{(y_j^{t+1})^2}{\theta_j} + \left(r\beta - \frac 32\right)\sum_{j=1}^n \log\left(\frac{\theta_j}{\vartheta_j}\right) - \sum_{j=1}^n \left(\frac{\theta_j}{\vartheta_j}\right)^r
	\right\}.
	\] 
	
	The passage from the signal to its increments is more challenging in dimensions $d\geq 2$ where the one-to-one correspondence no longer holds.  
	
	\begin{figure}
		\centerline{
			\includegraphics[width=5cm]{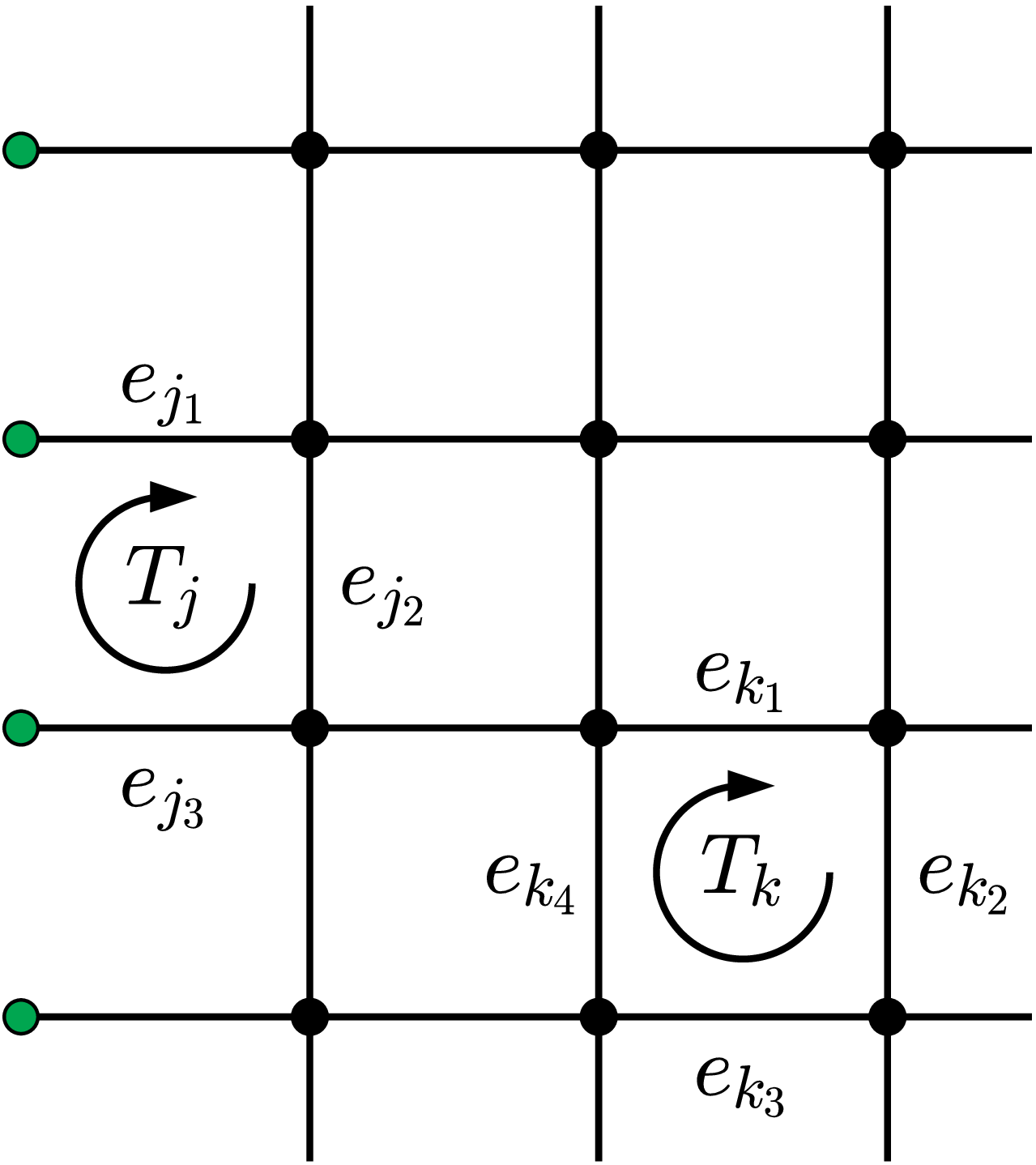}
		}
		\caption{\label{fig:circulation}Schematics of the circulation condition $\mM z=0$. In this figure, the black dots indicate the free nodes and the green dots are bound nodes in which the grid function is assumed to vanish. Between free nodes, no edge is defined, thus corresponding to a zero contribution to the circulation. If $z$ is a jump vector corresponding to a grid function $x$, the sums around edges of each loop must vanish.}
	\end{figure}
	
	To illustrate how to proceed, consider a quadrilateral graph, the nodes representing, e.g.,  the pixels in an image that we want to estimate,  with adjacent pixels connected by an edge, see Figure~\ref{fig:circulation}. Assume for simplicity that the values of the image vanish at the  boundary nodes, that we refer to as bound nodes, thus we are only interested in estimating the values at the remaining nodes, referred to as free nodes. Let $n_v$ be the number of the free nodes and $n_e$ the number of edges with at least one free node as an end point, referred to as free edges. 
	Let $\mL\in\R^{n_e\times n_v}$ denote the mapping from the free nodal values collected in the vector $x$ to the increments along free edges in the vector $y$,
	\begin{equation}\label{L}
	y = \mL x.
	\end{equation}
	Since the nodal values at the bound nodes, not included in the vector $x$, are equal to zero, the matrix $\mL$ has a trivial null space, i.e., ${\mathscr N}(\mL) = \{0\}$.
	
	Let $n_t$ denote the number of all loops $T_j$ in the graph, see Figure~\ref{fig:circulation}, including those defined by the edges between bound nodes,  and let $\mM\in\R^{n_e\times n_t}$ be the matrix computing the circulation around each loop by summing the increments over edges in clockwise order. If the increments along the edges correspond the nodal values, then the circulation in each element must vanish, i.e., $\mM y = 0$ or, equivalently,  $y \in{\mathscr N}(\mM)$ for every $y\in{\mathscr R}(\mL)$. %
	Since the edge increments associated with the nodal values are computed via the matrix $\mL$.
	Consequently, the matrices $\mL$ and $\mM$ define a short exact sequence,
	\[
	\{0\} \longrightarrow   \R^{n_v} \mathop{\longrightarrow}^{\mL} \R^{n_e}   \mathop{\longrightarrow}^{\mM} \R^{n_t} \longrightarrow  \{0\}. 
	\] 
	
	To define a prior promoting sparse increments, we consider a conditionally Gaussian prior model in terms of the increments $y_j$ along the free edges, written concisely as
	\begin{equation}\label{y prior}
	\pi_{y,\theta}(y,\theta) = 
	\pi_{y\mid\theta}(y\mid\theta)\pi_{\theta}(\theta)  \propto{\rm exp}\left(-\frac 12 \sum_{j=1}^{n_e} \frac {y_j^2}{\theta_j} + \phi(\theta)\right),
	\end{equation}
	where the function $\phi(\theta) = \phi_{r,\beta}(\theta)$ does not depend on $y$.
	It follows from the definition of the increments in terms of the nodal values (\ref{L}) that $y\in{\mathscr R}(\mL) = {\mathscr N}(\mM)$, therefore the support of the prior is restricted to ${\mathscr N}(\mM)$. 
	Introducing the auxiliary variable
	\[
	\beta = \mD^{-1/2}_\theta y,
	\]
	where $\mD_\theta\in\R^{n_e\times n_e}$ is the diagonal matrix with entries $\theta_j$, the compatibility condition on $y$ can be written in terms of $\beta$ as
	\begin{equation}\label{L theta}
	\beta \in {\mathscr R}(\mL_\theta),\quad  \mL_\theta =  \mD^{-1/2}_\theta \mL.
	\end{equation}
	Consider the QR factorization of $\mL_\theta$,
	\[
	\mL_\theta = \mQ\mR = \left[\begin{array}{cc} \mQ_1 & \mQ_2\end{array}\right]\left[\begin{array}{c} \mR_1 \\ \mO\end{array}\right],
	\]
	where the orthogonal matrix $\mQ$ is partitioned in the two blocks, $\mQ_1 \in\R^{n_e\times n_v}$, $\mQ_2 \in\R^{n_e\times (n_e-n_v)}$, and $\mR_1\in\R^{n_v\times n_v}$ is upper triangular and nonsingular because $\mL$ is of full rank, and $\mO$ is the zero matrix of size$(n_e-n_v)\times n_v$. For any $ \beta \in {\mathscr R}(\mL_\theta)$, there exists $x \in \R^{n_v}$ such that
	\begin{equation}\label{beta x}
	\beta = \mL_\theta x = \mQ \mR x,
	\end{equation}
	hence, by multiplying both sides of (\ref{beta x}) by the transpose of $\mQ$, we get 
	\[
	Q^\mT\beta = \left[\begin{array}{c} \mQ_1^\mT\beta \\ \mQ_2^\mT\beta\end{array}\right] = \left[\begin{array}{c} \mR_1 x \\ 0\end{array}\right],
	\] 
	or, equivalently,
	\[
	\mR_1^{-1} \mQ_1^\mT\beta = x,\quad \mQ_2^\mT \beta = 0.
	\]
	Therefore we can express the compatibility condition in terms of the auxiliary variable as 
	\[
	\beta\in{\mathscr N}(\mQ_2^\mT) = {\mathscr H}.
	\] 
	The posterior density for the prior (\ref{y prior}) on the increments, 
	\begin{equation}\label{posterior beta}
	\widetilde\pi_{\beta\mid b,\theta}(\beta\mid b,\theta) \propto {\rm exp}\left(-\frac{1}{2}\| b - \mA \mR_1^{-1}\mQ_1^\mT\beta\|^2  - \frac 12 \|\beta\|^2 +\phi(\theta)\right)
	\end{equation}
	if we neglect the compatibility conditions, when restricted to the subspace ${\mathscr H}$ becomes
	\begin{eqnarray}\label{posterior beta 2}
	\pi_{\beta\mid b,\theta}(\beta\mid b,\theta) &=&   \widetilde\pi_{\rm post}(\beta\mid b,\theta)\otimes \delta_{\mathscr H}(\beta) \\
	&\propto& {\rm exp}\left(-\frac{1}{2}\| b - \mA \mR_1^{-1}\mQ_1^\mT\beta\|^2  - \frac 12 \|\mQ_1^\mT \beta\|^2 +\phi(\theta)\right)\bigg|_{\mQ_2^\mT\beta = 0},\nonumber
	\end{eqnarray}
	where $\delta_{\mathscr H}$ is the singular measure concentrated on ${\mathscr H}$. 
	The following theorem shows that it is possible to carry out the iterations of the IAS algorithm for the posterior (\ref{posterior beta 2}) working directly with  (\ref{posterior beta}).
	
	\begin{theorem}
		The vector $\beta_*$ that maximizes (\ref{posterior beta}) satisfies $\mQ_2^\mT\beta_* = 0$, therefore also 
		maximizes (\ref{posterior beta 2}). Moreover, $\beta_*$ can be found by minimizing the expression
		\begin{equation}\label{pseudoinv}
		F(\beta) =  \frac{1}{2}\| b - \mA  \mL_\theta ^\dagger\beta\|^2  + \frac 12 \|\beta\|^2,
		\end{equation}
		where $\mL_\theta^\dagger$ is the pseudoniverse of $\mL_\theta$.
	\end{theorem}  
	
	{\em Proof.} From the observation that 
	\[
	\|\beta\|^2 = \| \mQ^\mT \beta \|^2 = \| \mQ_1^\mT \beta\|^2 +  \| \mQ_2^\mT \beta\|^2,
	\]
	it follows that 
	\begin{eqnarray*}
		&&\frac{1}{2}\| b - \mA \mR_1^{-1}\mQ_1^\mT\beta\|^2  + \frac 12 \|\beta\|^2 -\phi(\theta) \\
		&&\phantom{XXX}  = 
		\frac{1}{2}\| b - \mA \mR_1^{-1}\mQ_1^\mT\beta\|^2  + \frac 12 \|\mQ_1^\mT\beta\|^2 +
		\frac 12 \|\mQ_2^\mT\beta\|^2-\phi(\theta).
	\end{eqnarray*}
	and its minimizer, for fixed $\theta$, is 
	\[
	\beta_* = {\rm argmin}\left\{\frac{1}{2}\| b - \mA \mR_1^{-1}\mQ_1^\mT\beta\|^2  + \frac 12 \|\mQ_1^\mT\beta\|^2\right\},\quad
	\mQ_2^\mT\beta^* = 0,
	\]  
	which also maximizes (\ref{posterior beta 2}). Moreover, for $\beta^*$ such that $\mQ_2^\mT \beta^*=0$,
	\[
	\mR_1^{-1}\mQ_1^\mT\beta_* = \mL_\theta^\dagger\beta_*,
	\] 
	which completes the proof.$\quad\Box$ 
	
	The previous theorem shows that to find the MAP estimate, is not necessary to form the matrix $\mM$ or to compute the QR factorization of the matrix $\mL_\theta$.
	Instead, it suffices to solve the linear system
	\begin{equation}\label{beta}
	\left[\begin{array}{c} \mA \mL_\theta^\dagger \\ \mI\end{array}\right]  \beta= \left[\begin{array}{c} b \\ 0\end{array}\right]
	\end{equation}
	in the least squares sense, because its solution automatically satisfies $\mQ_2^\mT\beta = 0$, thus guaranteeing the existence of a vector $x$ such that  (\ref{beta x}) holds. The vector $y = \mD^{1/2}_\theta \beta$ satisfies the compatibility condition $\mM y=0$, representing feasible and consistent increments along the edges. 
	When resorting to the RKS approximation of the update of the signal inside the IAS iteration, we need to have a procedure to multiply a vector $\beta$ by the matrix  $ \mA \mL_\theta^\dagger$ and its transpose.
	The matrix-vector product of $\beta$ with $ \mA \mL_\theta^\dagger$  can be computed by first solving $\mL_\theta \alpha = \beta$ for $\alpha$ in the least squares sense, then multiplying $\alpha$ by $\mA$. 
	To evaluate the product of a vector $z$  we observe that the transpose of $ \mA \mL_\theta^\dagger$ is
	\[
	\big(  \mL_\theta^\dagger\big)^\mT \mA^\mT = \mL_\theta  \big(\mL_\theta^\mT \mL_\theta\big)^{-1}\mA^\mT,
	\]
	where, fortunately, in our case $\mL_\theta^\mT \mL_\theta$ is very sparse. Therefore we solve $\big(\mL_\theta^\mT \mL_\theta\big) w =  \mA^\mT z$, then multiply the solution $w$ by $\mL_\theta$.

	\section{Computed examples}
	\label{sec:ex}
	In our evaluation of the performance of the local and global hybrid IAS algorithms, we focus on the following questions: 
	
	{\em Stability of the convexity condition in local IAS:} To monitor how the components behave with respect to the local convexity region, we run the local hybrid IAS and monitor the behavior of the index set $I$  in Algorithm \ref{alg:1}. In particular, we track the indices $j\in I$, pointing to components $x_j$ that enter the local convexity region, satisfying \eqref{bound hyb}, and check whether or not they remain in $I$ without forcing the bound constraint for $x_j$. 
	
	{\em Identification of the support:} It is of particular interest to see whether the local hybrid method identifies correctly the support of a generative signal, avoiding stopping at a local minimum that may miss some of the components in the support, as sometimes happens when the non-convex prior models are used.
	Likewise, with the global hybrid algorithm, we monitor the indices corresponding to the variances in the convexity region at the switching iteration $\overline {t}$ and at the final iteration of Algorithm \ref{alg:2}.
	
	In our examples, the hyperpriors for the hybrid schemes are the gamma ($r^{(1)} = 1$) and the inverse gamma ($r^{(2)} = -1$). The performance of local and global hybrid IAS algorithms are also compared with the plain IAS with either the gamma or the inverse gamma hyperprior. In the global hybrid IAS algorithm, the switch to the non-convex model occurs at iteration $\overline{t}=10$.

	\paragraph{Example 1} The first test case is a one-dimensional deconvolution problem. The generative model is a piecewise constant signal  $f:[0,1]\to \R$,  $f(0)=0$, and the data consist of a few discrete noisy observations, 
	\[
	b_j = \int_0^1 A(s_j - t) f(t) dt + \varepsilon_j,\quad 1\leq j\leq m,  \quad A(t) = \left(\frac{J_1(\kappa|t|)}{\kappa|t|}\right)^2, 
	\]
	where $J_1$ is the Bessel function of the first kind and $\kappa$ is a scalar controlling the width of the kernel that we set to $\kappa = 40$, yielding a kernel with ${\rm FWHM} = 0.08$. We discretize the integral as
	\[
	\int_0^1 A(s_j - t) f(t) dt \approx  \sum_{j=1}^n w_k  A(s_j-t_k)f(t_k), \quad 1\leq k\leq n,
	\] 
	where $t_k = (k-1)/(n-1)$ and the $w_k$'s are the trapezoidal quadrature weights.
	We generate the data with a dense discretization with $n=n_{\rm dense} = 1253$, while in the forward model used for solving the inverse problem, we set $n=500$. The observation points are given by $s_j = (4+ j)/100$,  $1\leq j\leq m = 91$, and  the additive noise is assumed to be scaled white noise, with standard deviation $\sigma$ set to  2\% of the maximum of the noiseless generated signal. We denote $x_j = f(t_j)$.  Figure~\ref{fig:1D_Data} shows the generative signal and the data.
	
	While the generative signal, a piecewise constant function, is not sparse, it admits a sparse representation in terms of its increments  $z_j = x_j - x_{j-1}$ over the interval of definition. Assuming  $x_0=0$, then 
	\begin{equation}\label{difference}
	z = \mL  x\,,\quad \mL=\left[\begin{array}{cccc}
	1&0&\ldots&0\\
	-1&1&\ldots&0\\
	&&\ddots&\\
	0&\ldots&-1&1\\
	\end{array}\right]
	\in \R^{n\times n },
	\end{equation}
	hence
	\begin{equation}
	x = \mC z\quad\mathrm{with}\quad  \mC = \mL^{-1} = \left[\begin{array}{cccc}
	1&0&\ldots&0\\
	1&1&\ldots&0\\
	\vdots &&\ddots&\\
	1&\ldots&1&1\\
	\end{array}\right]\in\R^{n\times n}.
	\end{equation}
	Therefore our inverse problem is to estimate the vector $z$, assumed to be sparse, from the data vector $b$, given the forward model 
	\begin{equation}
	b = \mA\mC z + e, \quad \varepsilon \sim \mathcal{N}(0,\sigma^2 \mI), \quad \mA_{jk} = w_k A(s_j - t_k).
	\end{equation}
	
	\begin{figure}
		\centerline{\includegraphics[width=4.8cm]{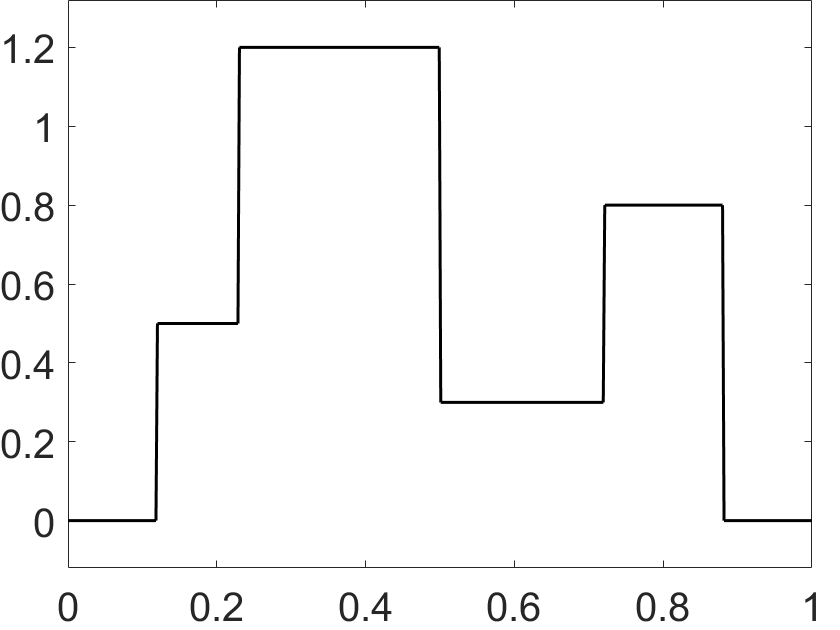} \quad \includegraphics[width=4.8cm]{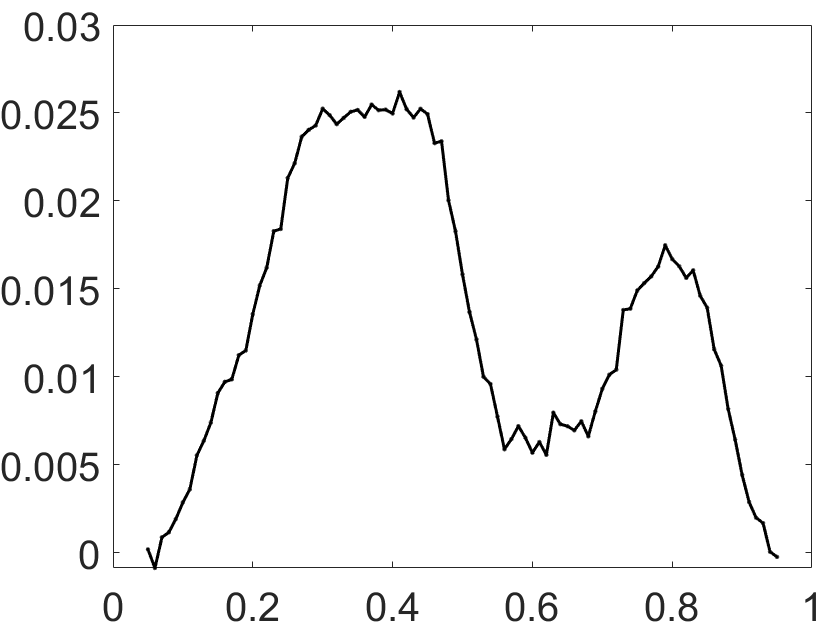}
		}
		\caption{\label{fig:1D_Data} Left: The generative model. Right: the blurred and noisy data vector $b \in \R^{91}$.}
	\end{figure}
	
	The reconstruction results, together with the final variance vector $\theta$ and the number of CGLS steps per IAS iteration, are shown in Figure~\ref{fig:1D_rest}. 
	The locations of the first two increments in the generative signal are not easy to detect from the data, see Figure~\ref{fig:1D_Data}, and they are not sharply restored by the IAS algorithm with gamma hyperprior (see the first row of Figure~\ref{fig:1D_rest}), while the IAS with inverse gamma (second row) hyperprior lumps the increments, stopping at a local minimizer that corresponds to a simpler profile.
	
	The reconstruction with the local hybrid algorithm is shown in the first panel of the third row of Figure~\ref{fig:1D_rest}. The middle panel of the same row shows in dashed blue
	the components of $\theta$ that follow the inverse gamma distribution at the last iteration of IAS, and in solid red those that never switch from the gamma distribution. The effect of changing to the inverse gamma is a cleaner background. The global hybrid hyperprior returns a sharp restoration of the signal, as shown in the first panel of the fourth row of Figure~\ref{fig:1D_rest}, with the five jumps accurately identified in the correct positions. In both cases, after a few steps, the number of CGLS steps per IAS iteration equals the number of increments detected, indicating that both hybrid IAS algorithms can determine very accurately the cardinality of the support: see also \cite{CPrSS}.

	\begin{figure}
		\centerline{\includegraphics[width=4cm]{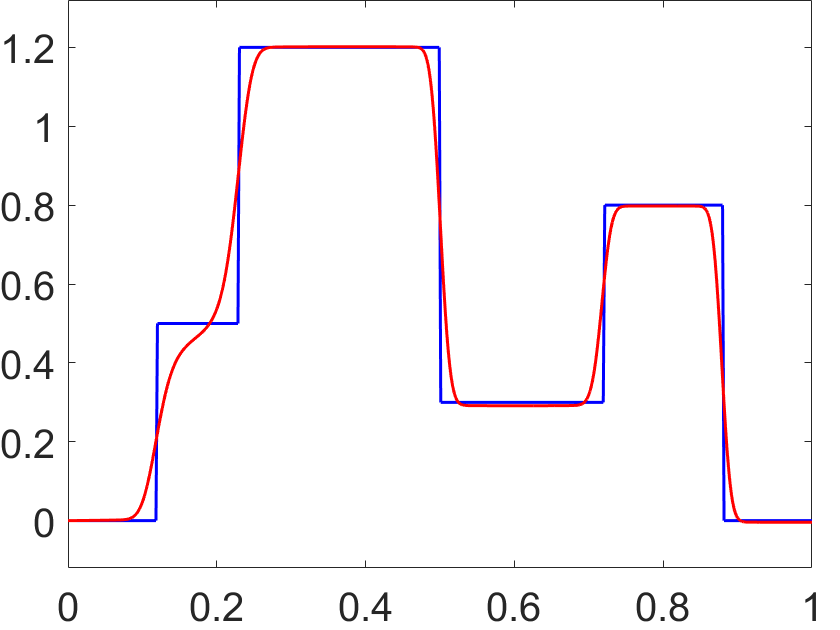} \quad \includegraphics[width=4cm]{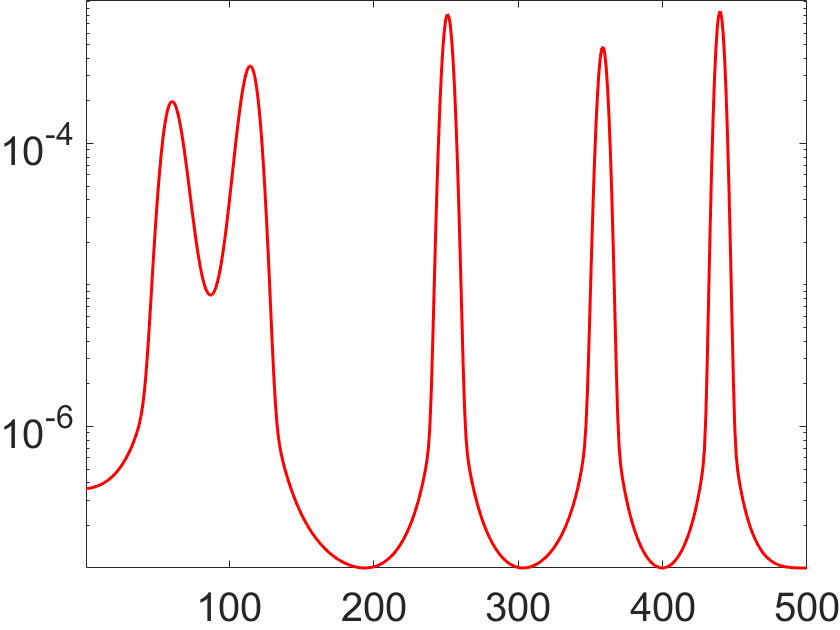}
			\quad \includegraphics[width=4cm]{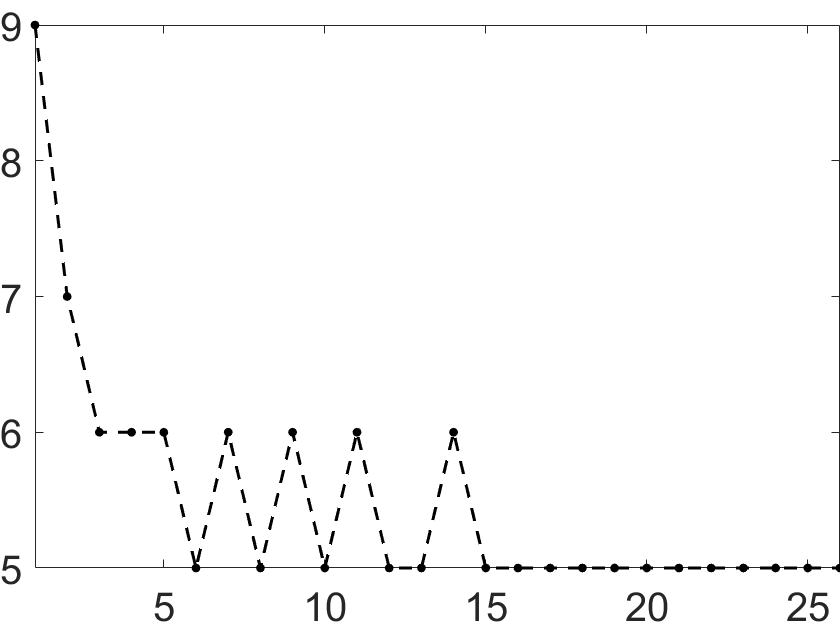}}
		\centerline{\includegraphics[width=4cm]{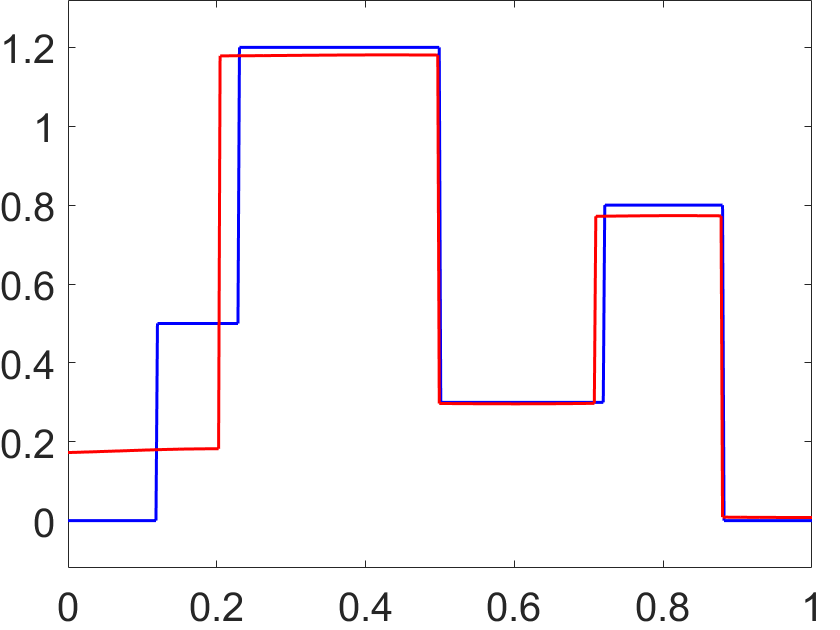} \quad \includegraphics[width=4cm]{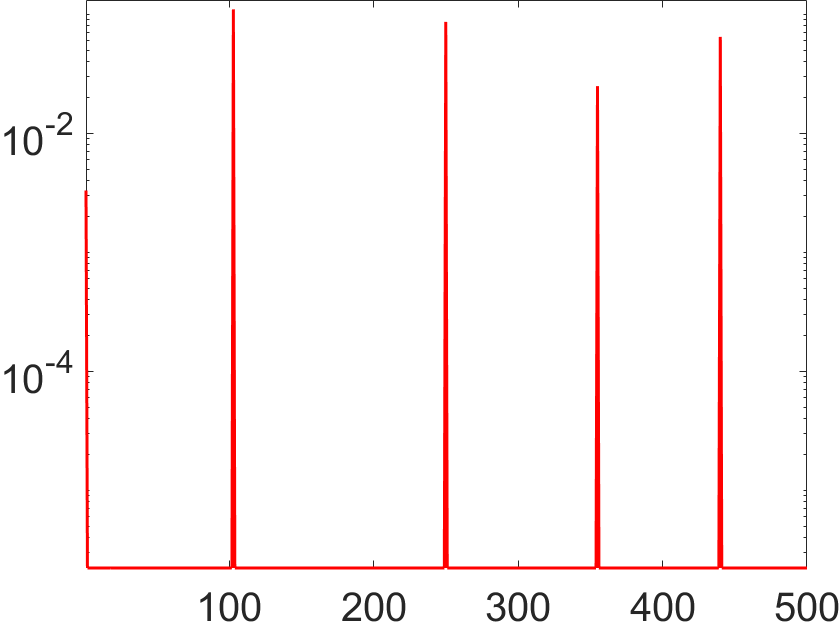}\quad \includegraphics[width=4cm]{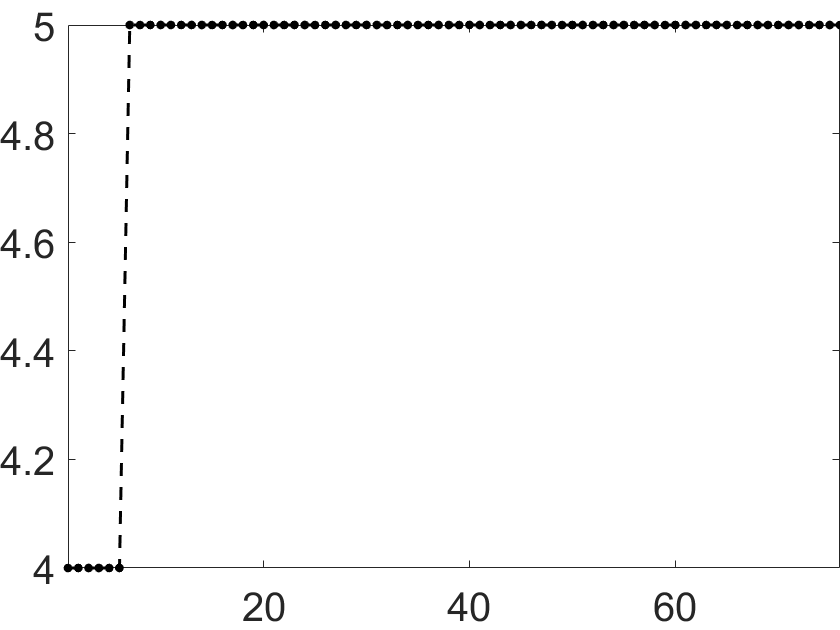}}
		\centerline{\includegraphics[width=4cm]{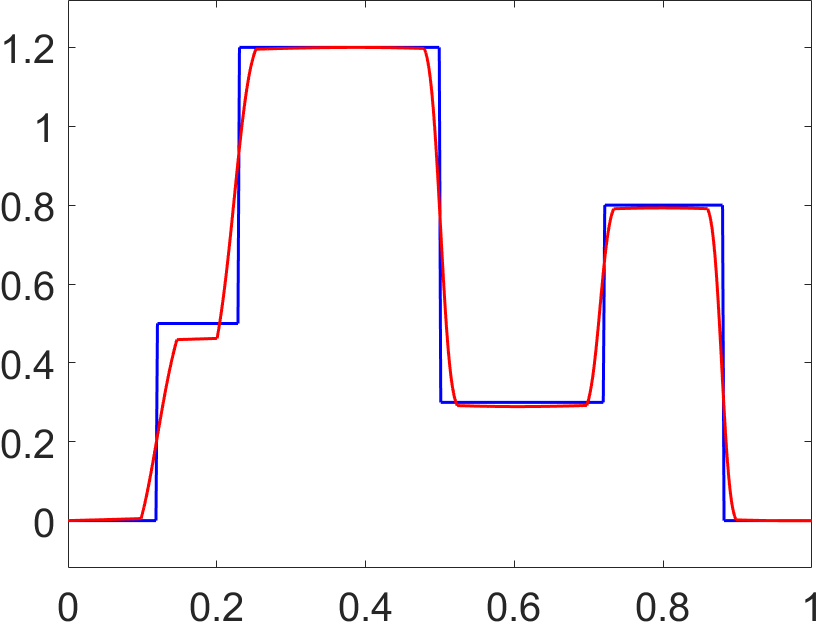} \quad \includegraphics[width=4cm]{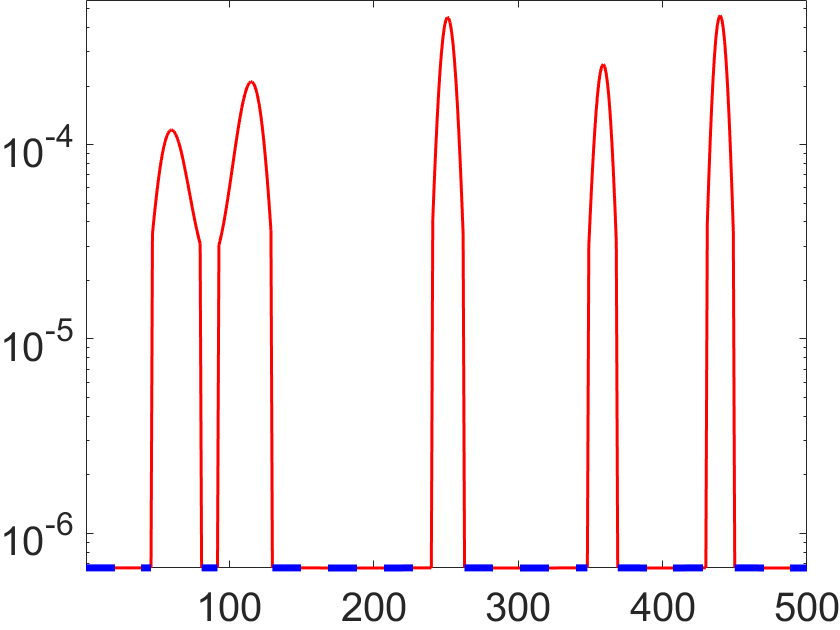}
			\quad \includegraphics[width=4cm]{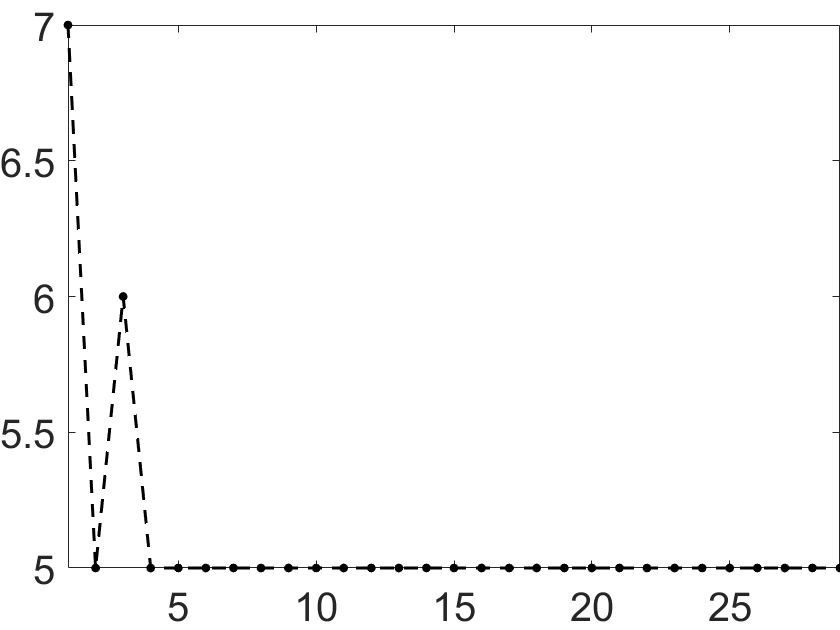}}
		\centerline{\includegraphics[width=4cm]{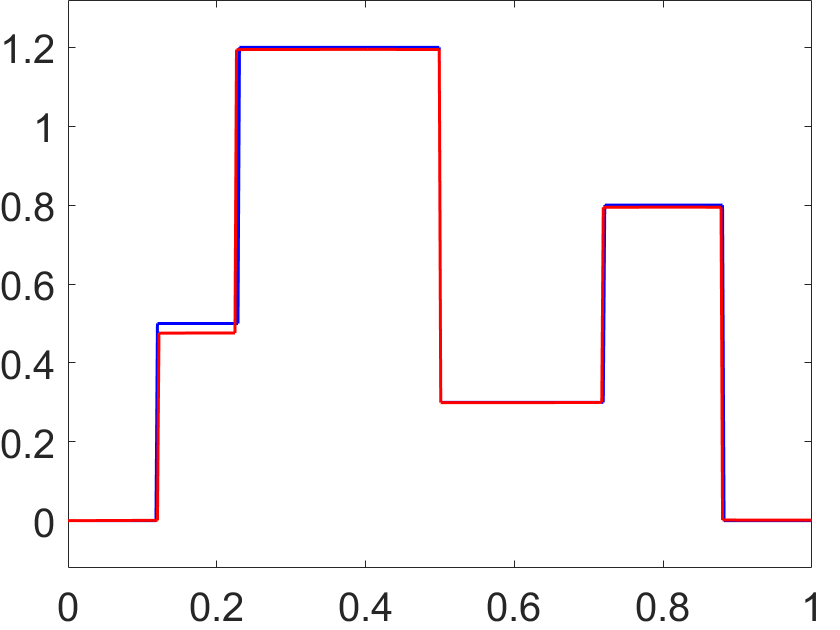} \quad \includegraphics[width=4cm]{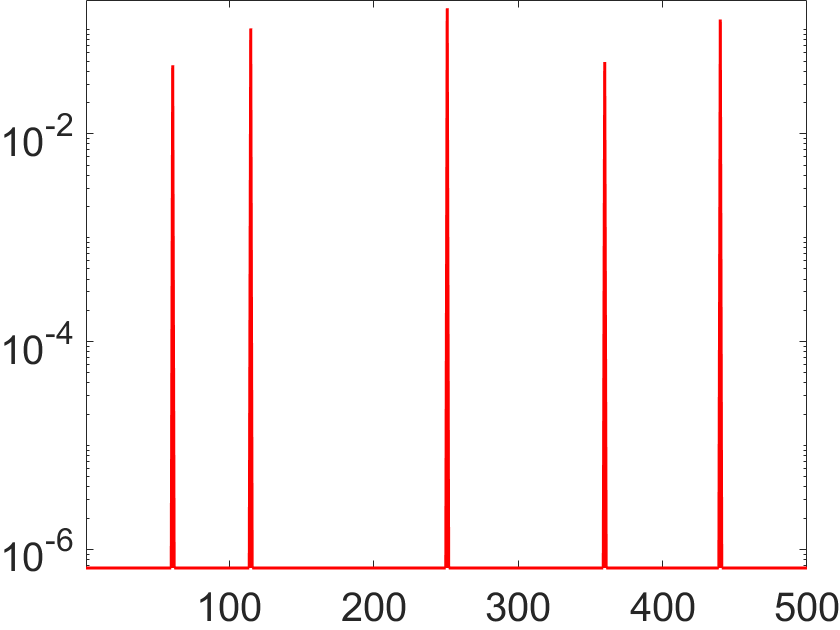}
			\quad \includegraphics[width=4cm]{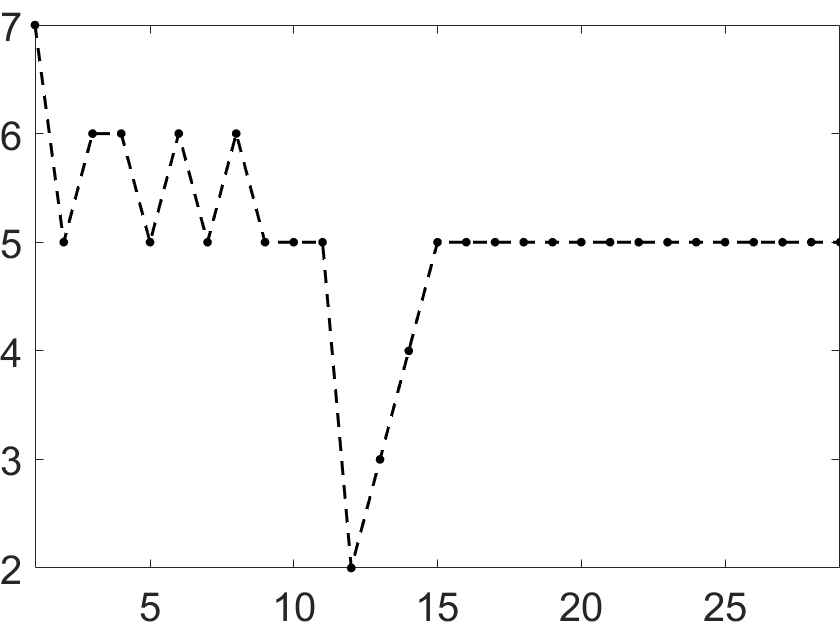}}
		\caption{\label{fig:1D_rest} Reconstruction of the signal via gamma, inverse gamma, local hybrid and global hybrid hyperprior (left), the hyperparameter $\theta$ (center) and the CGLS iterations per each IAS iteration (right). For the gamma hyperprior in the top row the parameter values are $\eta=10^{-2}$ and $\vartheta=10^{-5}$, for the inverse gamma hypeprior in the second row $\eta=-4.5$ and $\vartheta=10^{-5}$. The hybrid hyperpriors in the bottom rows inherit the parameters from the generative hyperpriors.}
	\end{figure}

	To address the stability of the convexity condition, we follow iteration by iteration the convexity condition, classifying each index in the sets $I$ (convexity condition satisfied) and its complement $I^c$ (condition not satisfied). The left panel of Figure \ref{fig:1D_vars}, where the indices in $I$ are marked in green, and those in $I^c$ in yellow, indicate that the set $I$ is monotonously increasing, that is, once a component enters the convexity region, it does not leave it, thus effectively removing the need for imposing the bound constraint \eqref{bound hyb}.
	
	The middle panel of Figure~\ref{fig:1D_vars} shows the variances $\theta_j$ in the global hybrid algorithm at the end of the iteration $\overline t - 1 = 9$, prior to switching to the inverse gamma model.  The components  $\theta_j$ satisfying the convexity bound, in dashed blue, are those for which the switch to the inverse gamma distribution does not compromise the convexity of the objective function. The panel on the right indicates for each component at each global hybrid IAS iteration whether it satisfies (green) or not (yellow) the convexity bound. Although in this case, unlike for the local hybrid IAS algorithm, the index set $I$ is not monotonically increasing, eventually the support is correctly detected to high accuracy.
	
	\begin{figure}
		\centerline{\includegraphics[width=4cm]{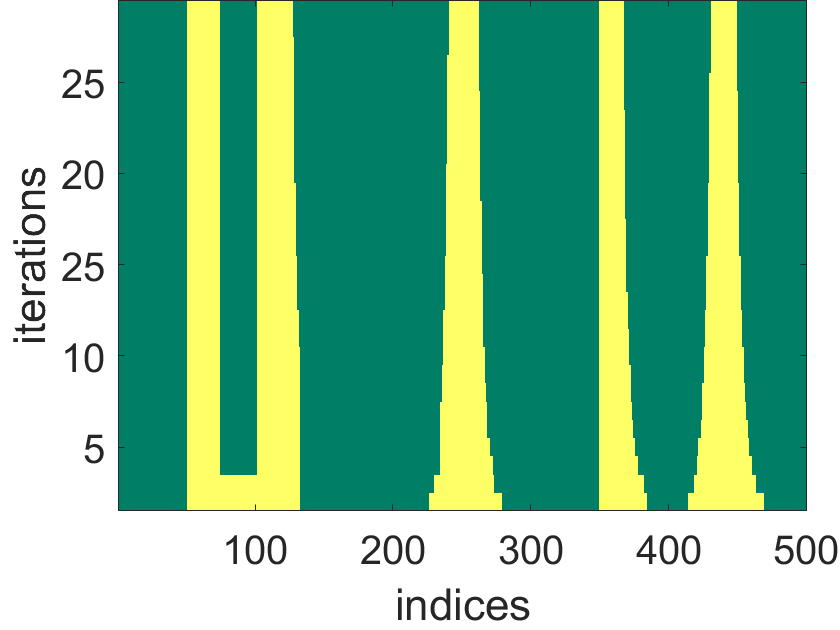} \quad \includegraphics[width=4cm]{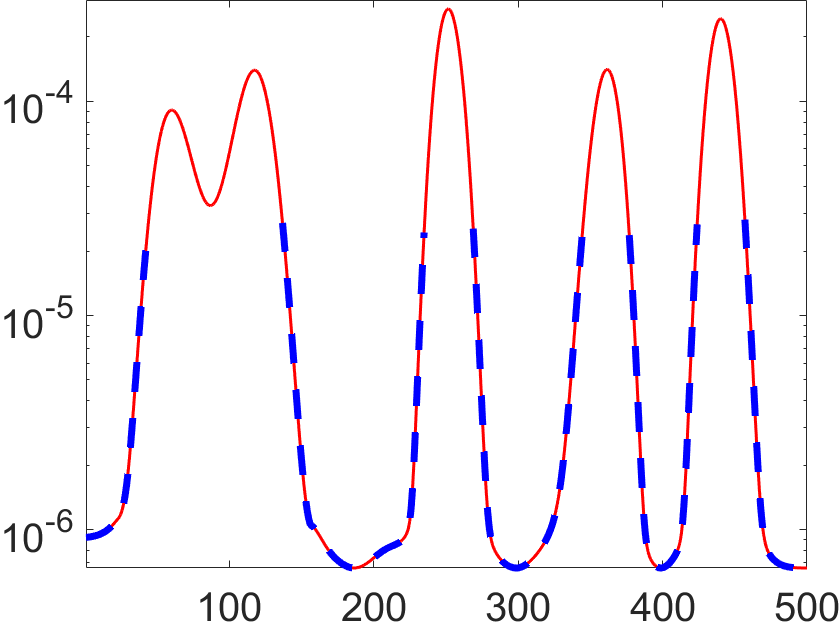}\quad \includegraphics[width=4cm]{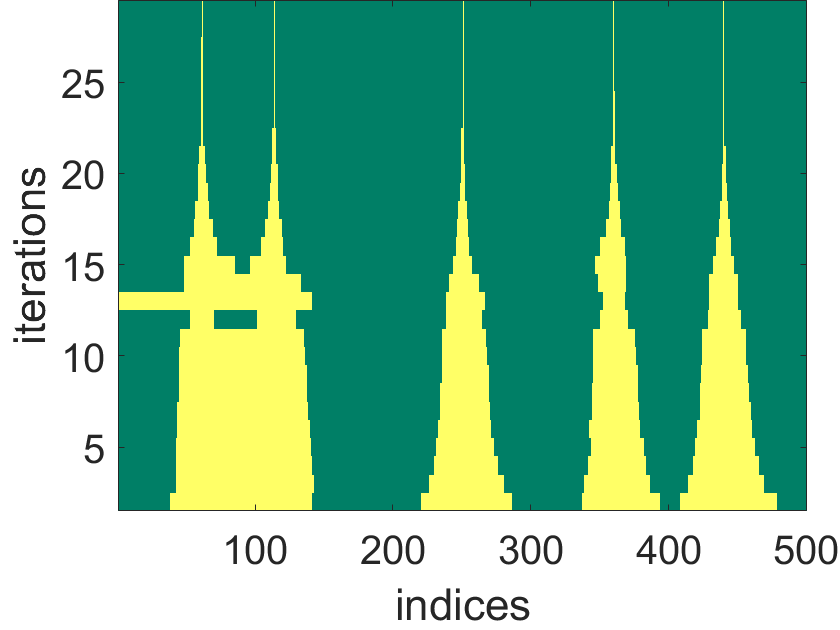}
		}
		\caption{\label{fig:1D_vars} Left: Pseudocolor image of the indices iteration by iteration of the local hybrid algorithm, green indicating the indices of those components that are in the convexity domain (index set $I$), and yellow those that are outside of it (index set $I^c$). Observe that when moving up, the yellow intervals shrink and the green ones increase, indicating stable convexity without the need to force the bound constraint \eqref{bound hyb}. Center: Variables $\theta_j$ in the global hybrid algorithm at the iteration $\overline t -1 = 9$, before the switch to inverse gamma model. The values in the convexity region are indicated in blue, the rest in red. Right: Pseudocolor image of the indices iteration by iteration in the global hybrid algorithm, green indicating the indices with components in the stability region. Observe that while the algorithm converges, correctly identifying the support, the index sets are not monotonous. In particular, after the switch ($\overline t = 10$), the discontinuities close to the left end of the interval create some confusion.}	
	\end{figure}

	\paragraph{Example 2} The second test case is an image restoration problem. Let $\Omega$ be a square compact region in $\R^2$  and $x$ be the generative image defined over $\Omega$. The discrete and noisy data consists of observation at points $q_j\in\Omega$ of a convolved version of the image,
	\begin{equation}
	\label{eq:obs data}
	b_j = \int_{\Omega}A(q_j,p')x(p') dp' + \varepsilon_j.
	\end{equation}
	with a  Gaussian convolution kernel 
	\begin{equation}
	\label{eq:gauss kern}
	A(p,p') = \frac{1}{2\pi w^2} \exp\left(-\frac{\|p-p'\|_2^2}{2 w^2}\right)\,,\quad\mathrm{with}\quad w= 0.015.
	\end{equation}
	The integral is discretized 
	over an $n \times n$ pixel grid with $n=136$, whereas the number of observation points is $m = 68\times 68$. The noiseless signal is corrupted by additive scaled white noise with standard deviation approximately $2\%$ of the maximum noiseless signal. We assume a priori sparsity of  the horizontal and vertical increments of the discrete image $x$, and implement the sparsity prior in the IAS algorithm according to the procedure detailed in Section \ref{sec:ias incr}. The original image, the observed data,  vertical and horizontal increments of the original image are shown in Figure \ref{fig:2D_data}.
	
	\begin{figure}
		\centerline{\includegraphics[width=2.9cm]{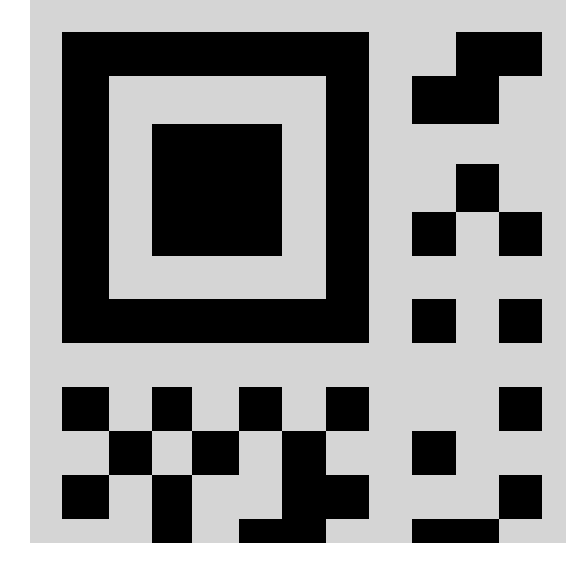} \quad \includegraphics[width=2.9cm]{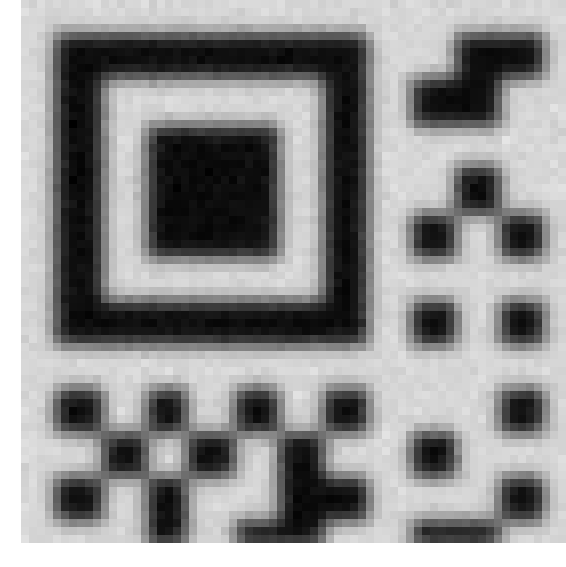} \quad \includegraphics[width=2.9cm]{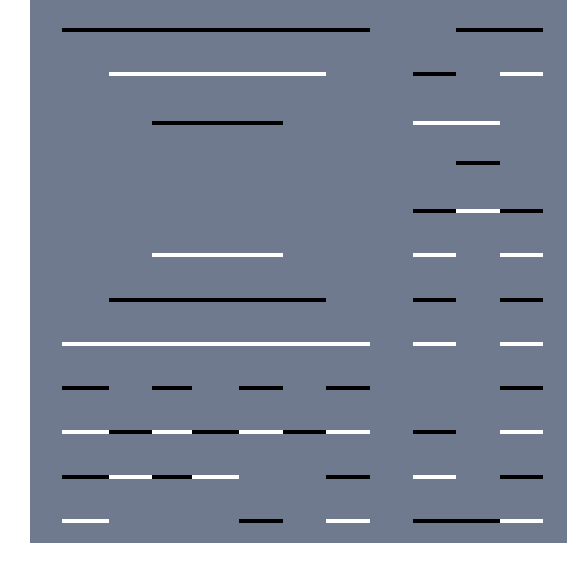} \quad \includegraphics[width=2.9cm]{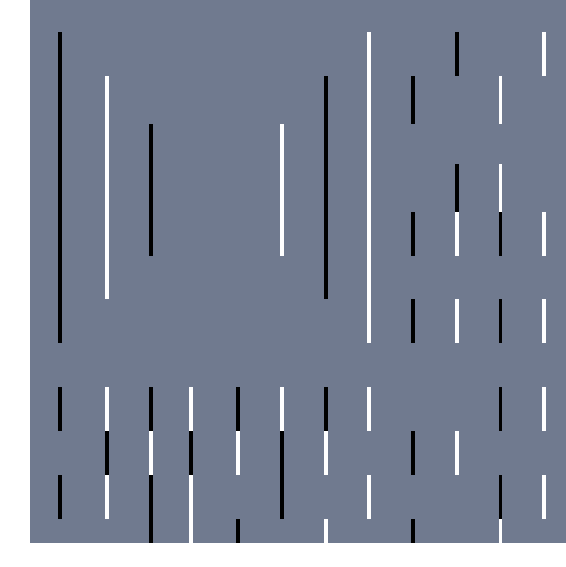}
		}
		\caption{\label{fig:2D_data} From left to right: original test image $x\in\R^{136 \times 136}$; observed data $b\in\R^{68 \times 68}$ corrupted by Gaussian blur and additive Gaussian noise; vector of horizonal increments of the image; and vertical increments.}
	\end{figure}
	
	The IAS is performed by constraining the values $x_j$ in the interval $[0,1]$, $1 \leq j \leq n^2$. More details on constrained IAS are given in \cite{CPrSS}.\\
	
	\vspace{0.2cm}
	
	The restored images computed by the IAS algorithm with gamma, inverse gamma, and by the local and global IAS algorithms are shown in Figure \ref{fig:rest_2d}.  Figure \ref{fig:var_2d} shows the logarithmic plot of variances $\theta_j$, the profile of the restorations along the dotted horizontal cut lines indicated in the reconstructions of Figure \ref{fig:rest_2d}, and the profile of the original image. Not surprisingly, the restoration using gamma hyperprior shows slightly rounded corners, while the algorithm based on inverse gamma hyperprior produces some staircasing artifacts along the edges. Both  effects are mitigated in the restorations computed with the hybrid IAS algorithms. The reconstruction of the global hybrid IAS algorithm is of remarkably high quality.
	
	\begin{figure}
		\centerline{
			\includegraphics[width=3.2cm]{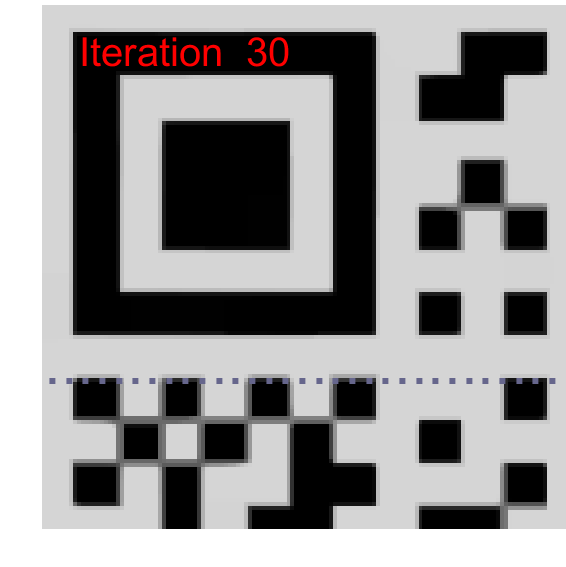} 
			\includegraphics[width=3.2cm]{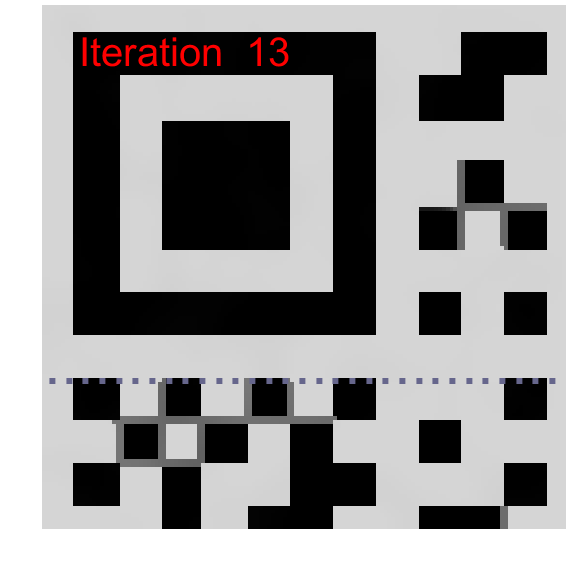}
			\includegraphics[width=3.2cm]{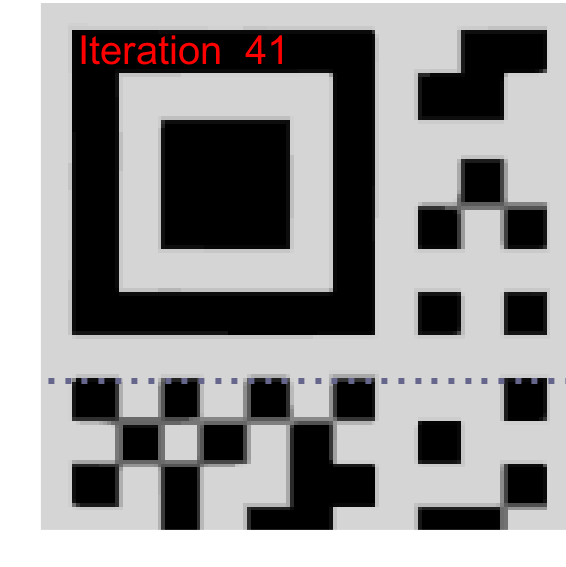}
			\includegraphics[width=3.2cm]{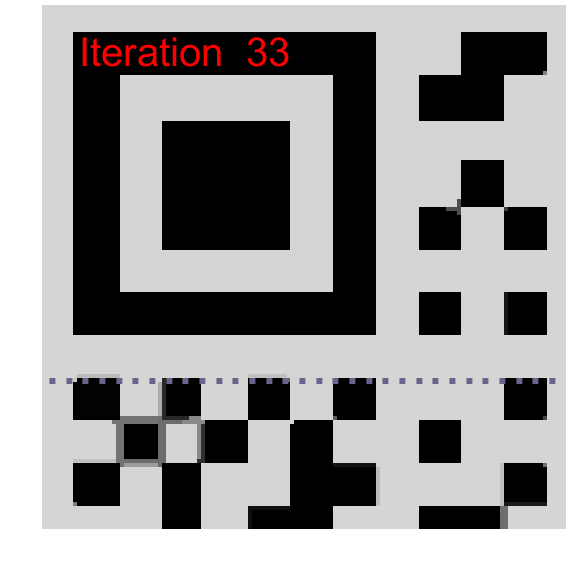}}
		\caption{ \label{fig:rest_2d}From left to right: Restored images by IAS algorithm based on gamma and inverse gamma hyperpriors, and by the local and global hybrid IAS algorithm using the combination of gamma and inverse gamma models. In the gamma hyperprior, the parameter values are $\eta=10^{-4}$ and $\vartheta=10^{-3}$, and in the inverse gamma hypeprior, $\eta=-6.5$ and $\vartheta=10^{-4}$. The hybrid hyperpriors  inherit these parameters from the generative hyperpriors. The dotted horizontal line indicates a cut across the reconstruction defining the profiles shwn in Figure~\ref{fig:var_2d}.}
		
	\end{figure}
	
	\begin{figure}
		\centerline{
			\includegraphics[width=3cm]{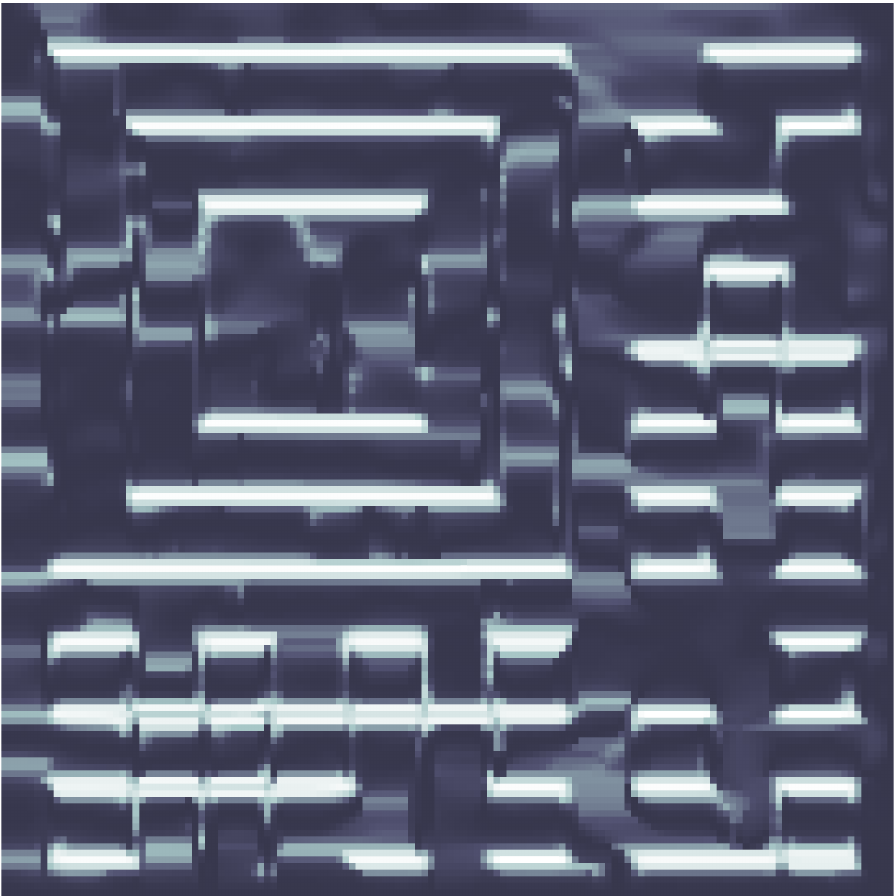}
			\quad\includegraphics[width=3.52cm]{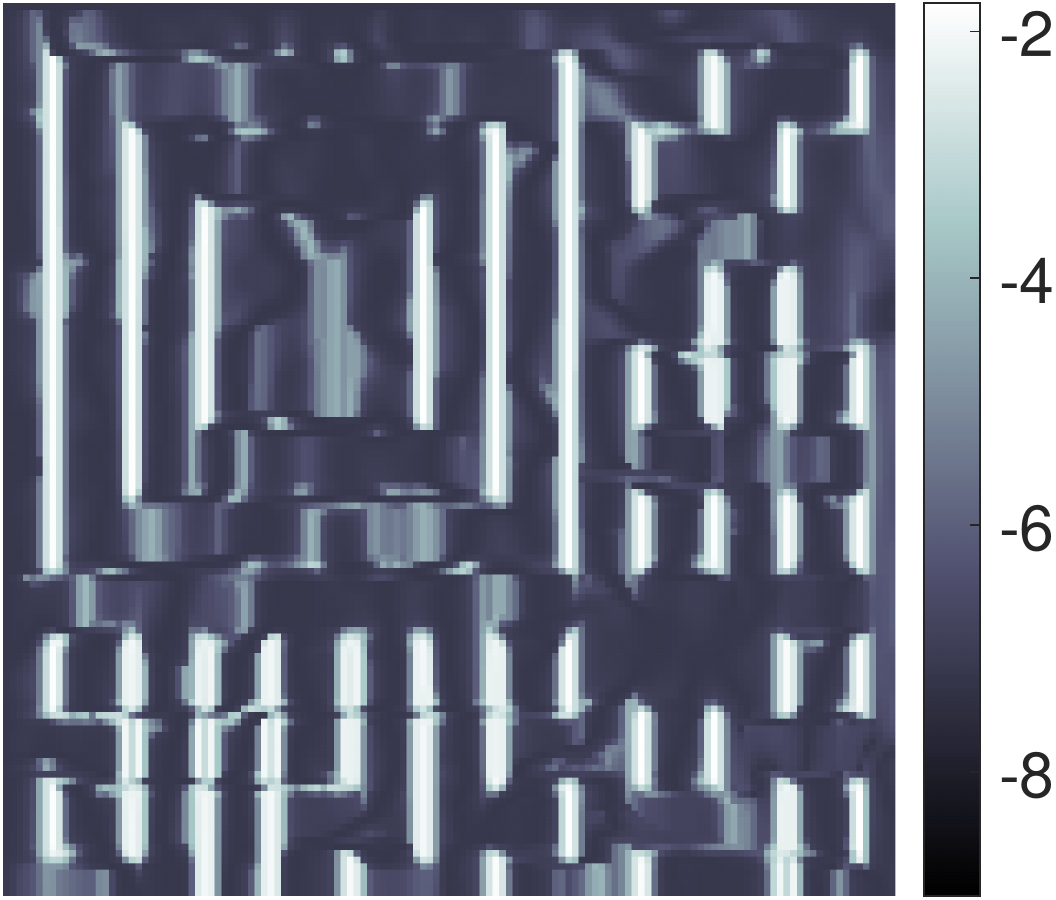} \quad 		\includegraphics[width=4cm]{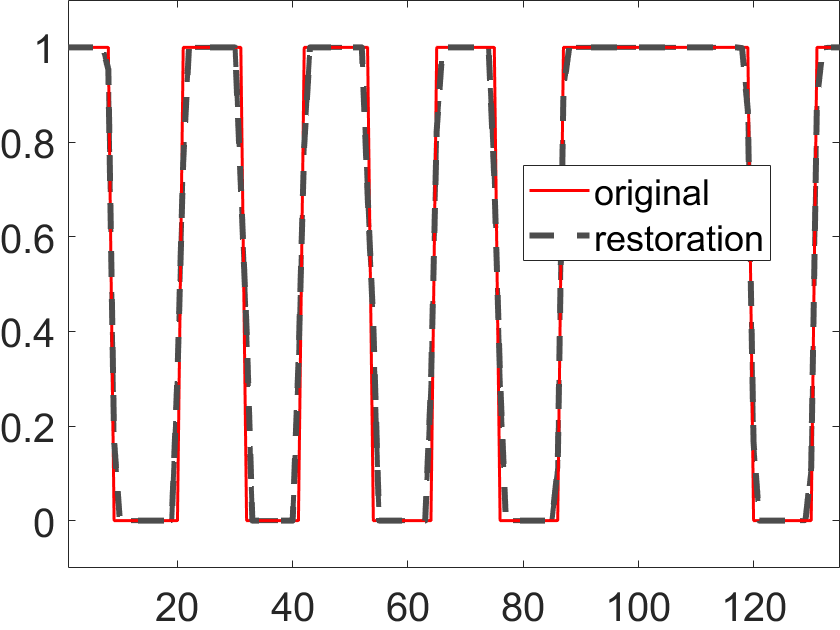}} 
		\centerline{
			\includegraphics[width=3cm]{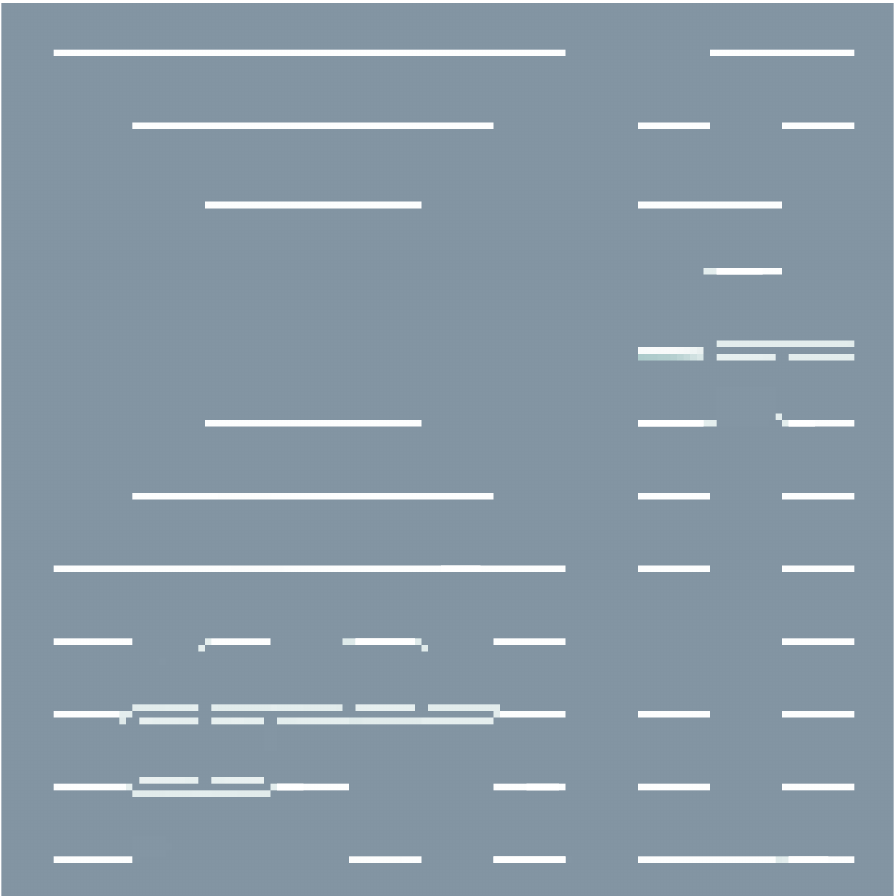}
			\quad\includegraphics[width=3.52cm]{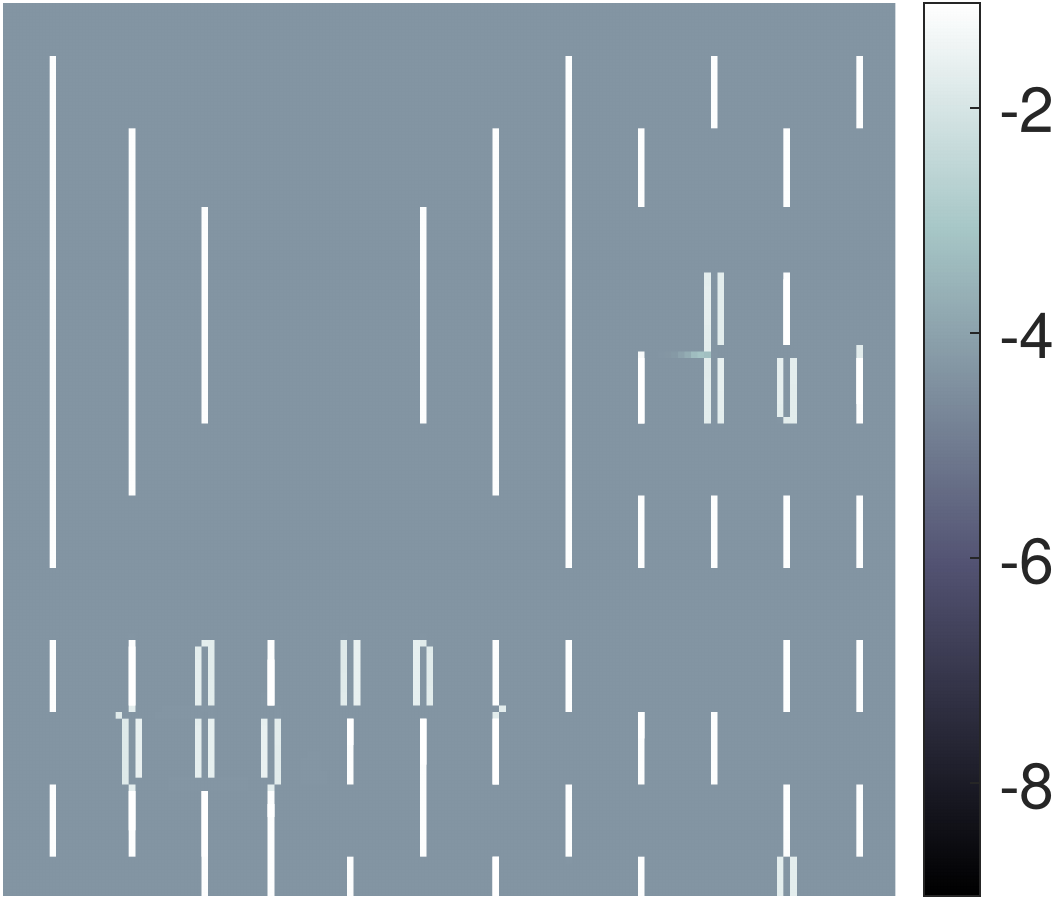} \quad 		\includegraphics[width=4cm]{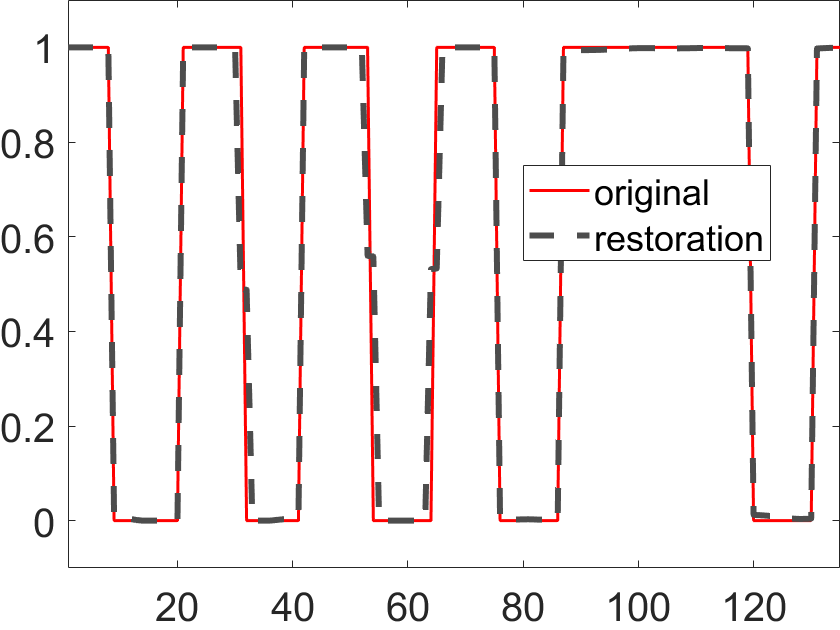}} 
		\centerline{
			\includegraphics[width=3cm]{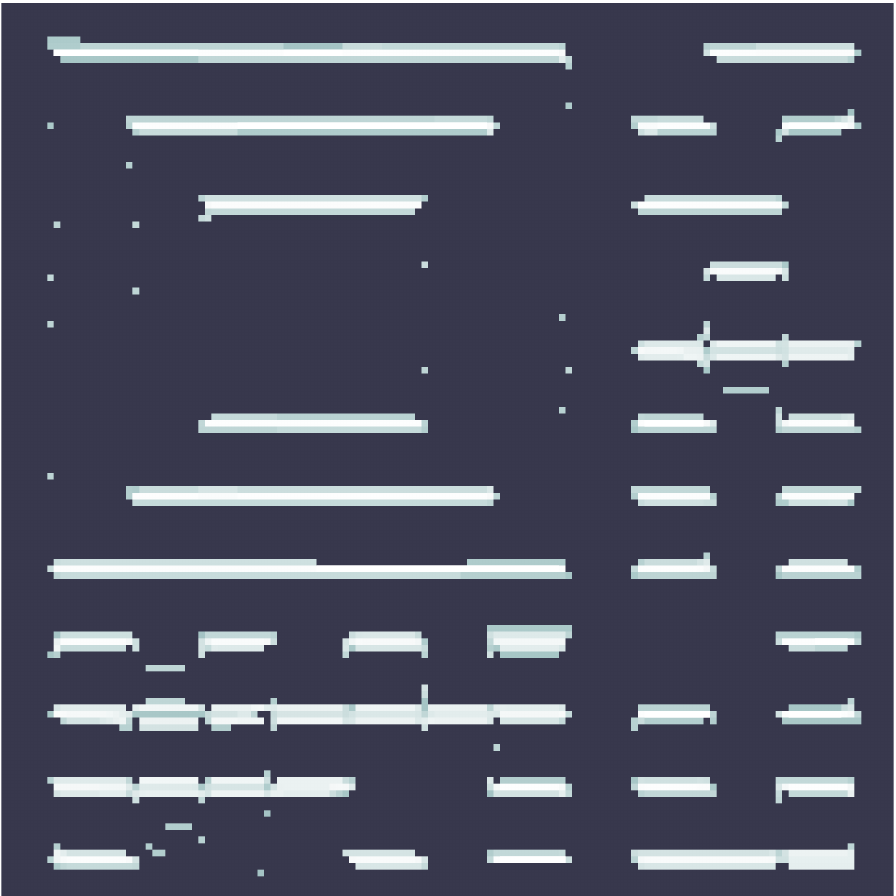}
			\quad\includegraphics[width=3.52cm]{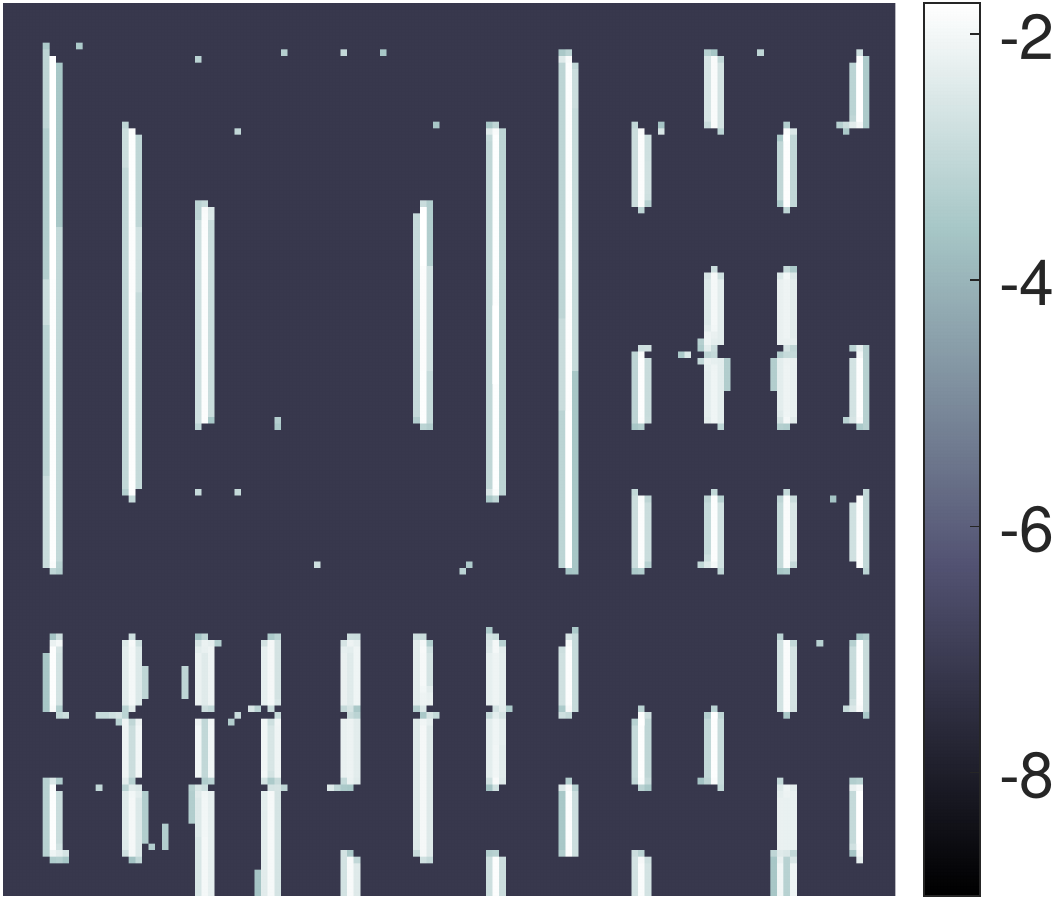} \quad 		\includegraphics[width=4cm]{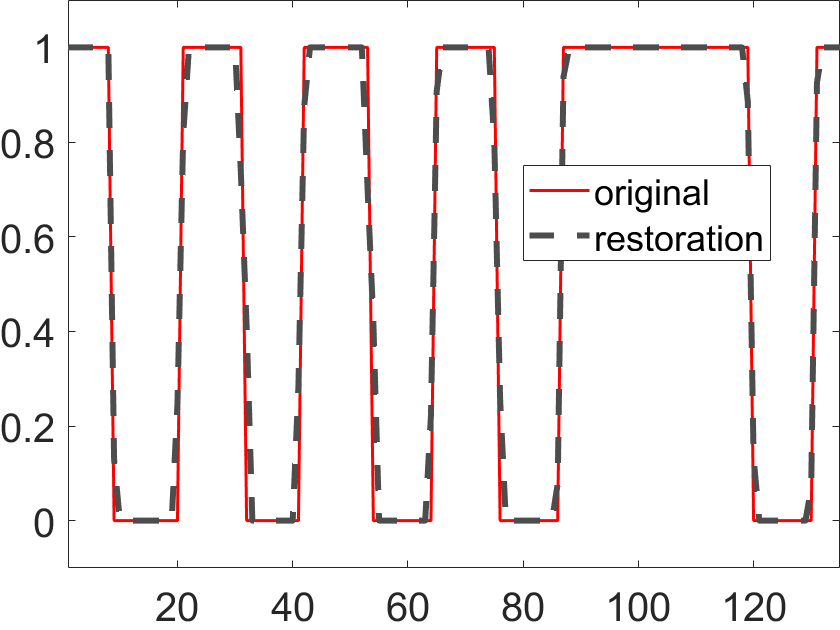}} 
		\centerline{
			\includegraphics[width=3cm]{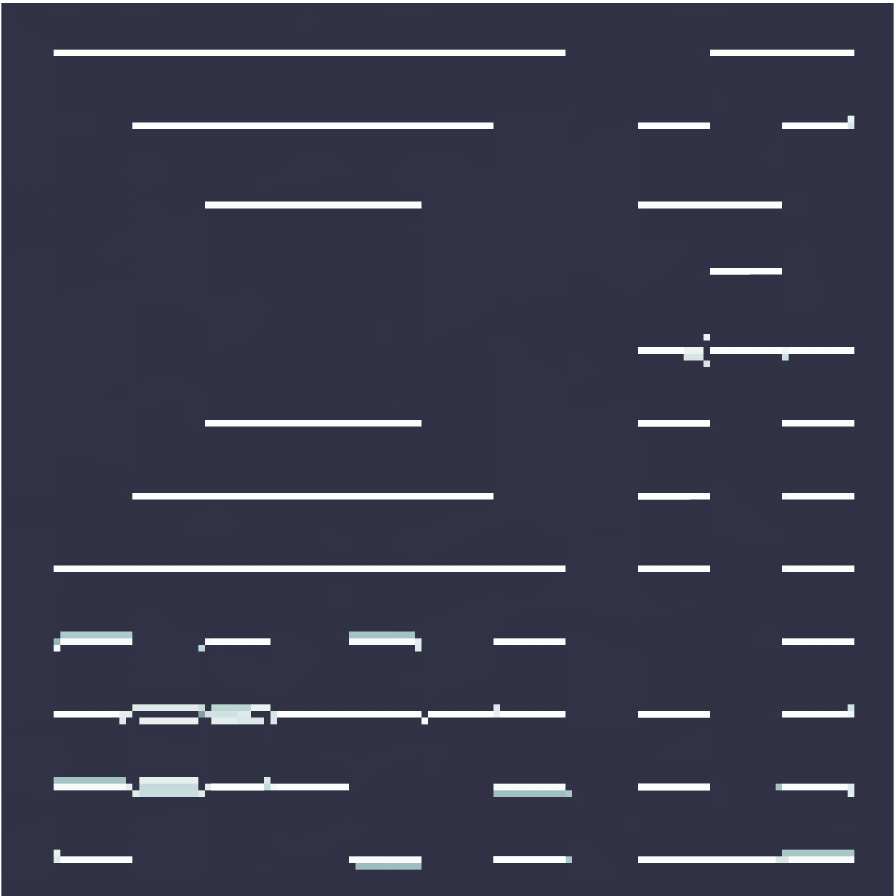}
			\quad\includegraphics[width=3.52cm]{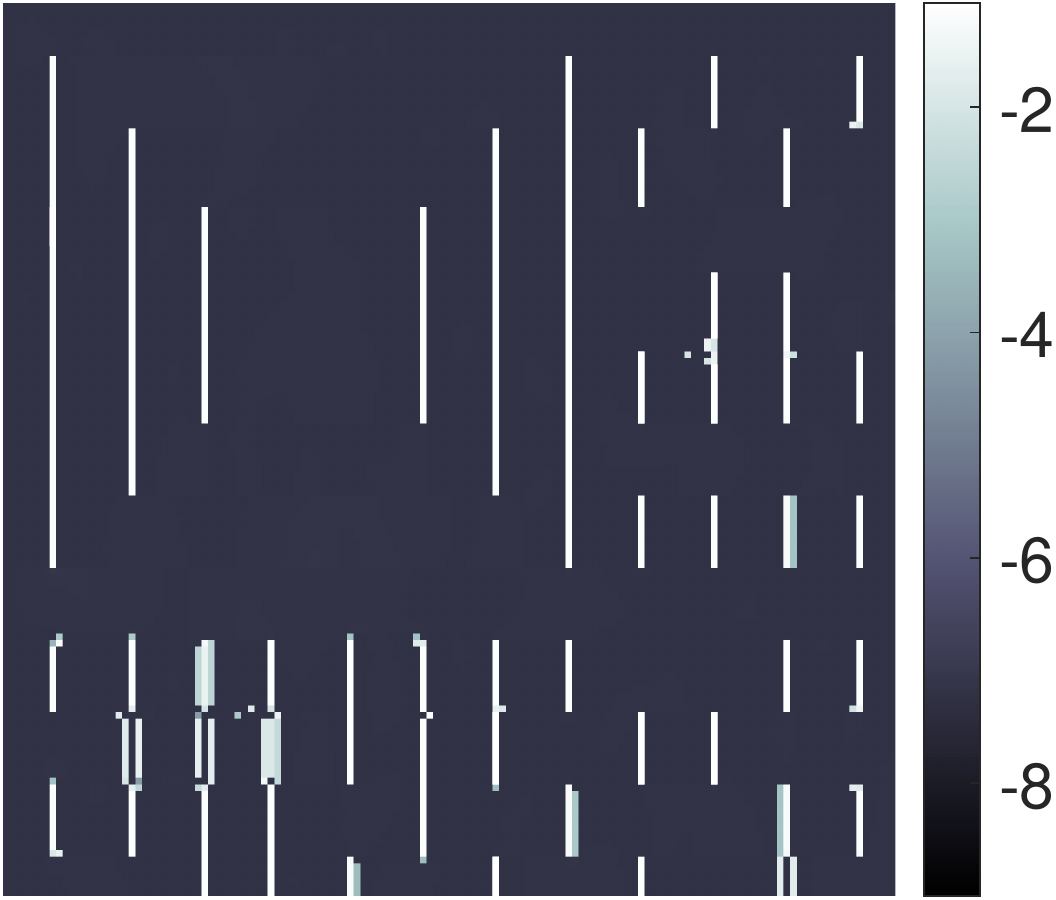} \quad 		\includegraphics[width=4cm]{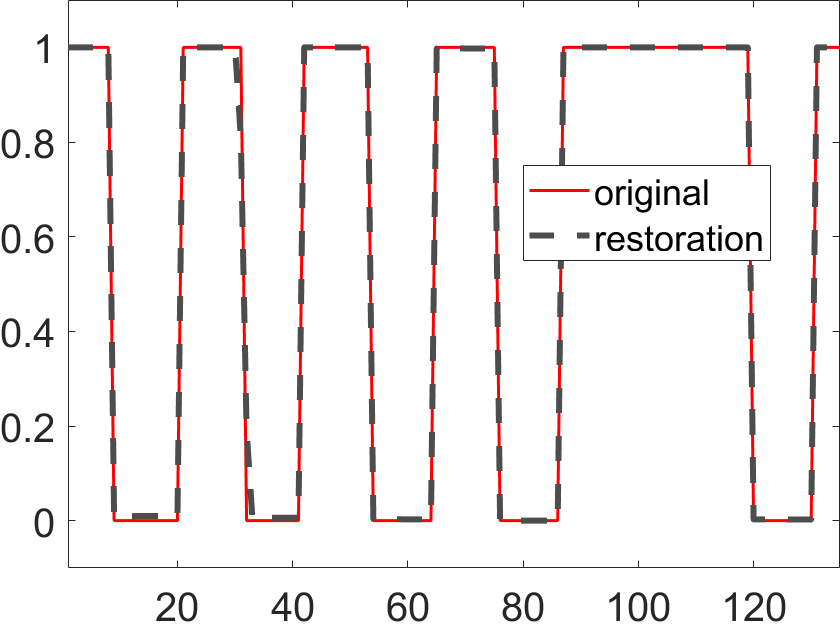}} 
		\caption{\label{fig:var_2d}From top to bottom: Logarithmic plots of variances corresponding to vertical and horizontal increments, and on the right, one-dimensional profiles extracted from the restorations in Figure \ref{fig:rest_2d} for the gamma, inverse gamma, local hybrid and global hybrid hyperprior.}
		
	\end{figure}
	
	The number of CGLS steps in each IAS iteration for the four models is reported in Figure \ref{fig:cgls qr}.
	
	\begin{figure}
		\centerline{		\includegraphics[width=5cm]{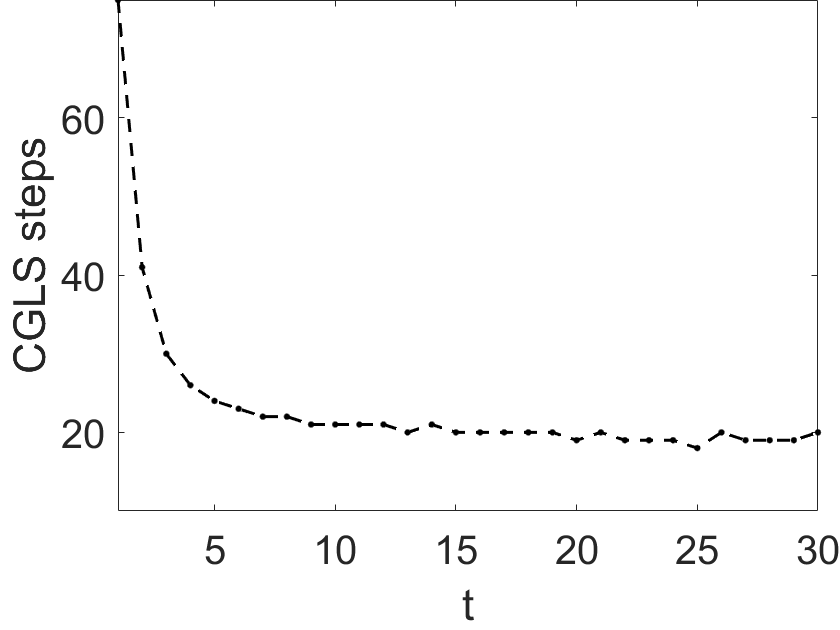}
			\quad
			\includegraphics[width=5cm]{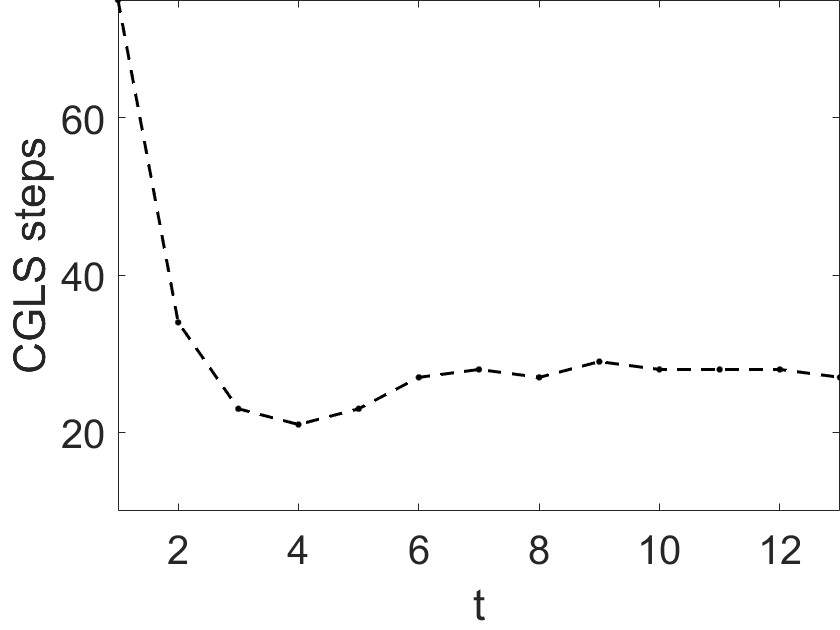}}
		
		\centerline{
			\includegraphics[width=5cm]{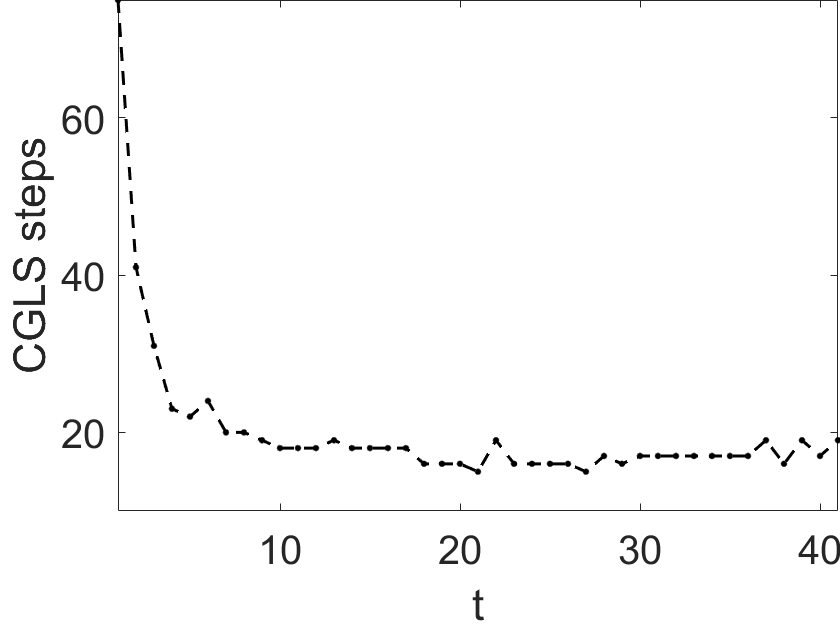}
			\quad 	\includegraphics[width=5cm]{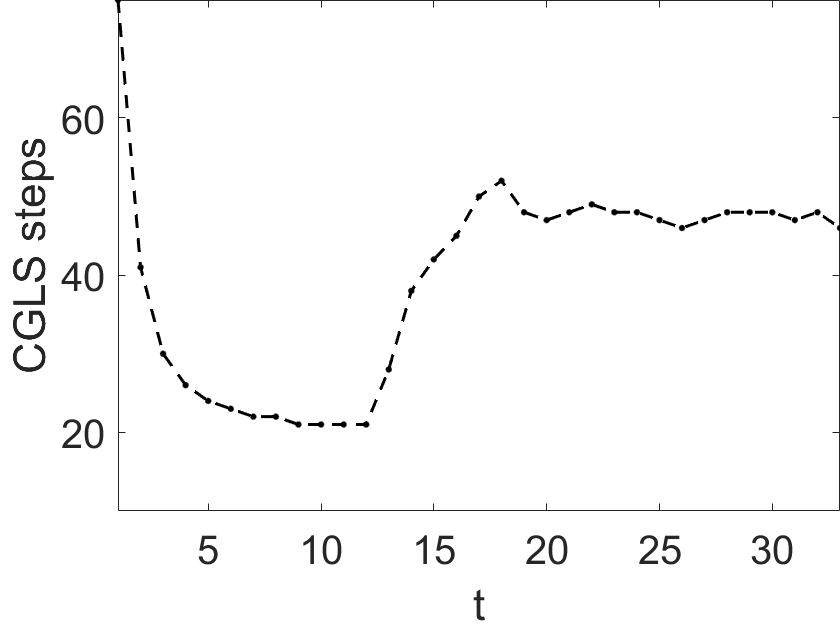}}
		\caption{Number of CGLS steps per outer iteration, from top to bottom and from left to right, for gamma, inverse gamma, local hybrid and global hybrid hyperprior.}
		\label{fig:cgls qr}
	\end{figure}

	The left panels of Figure \ref{fig:qr ind} display pseudocolor images of the indices of the variances $\theta_j$ of the horizontal and vertical increments at the last iteration of local hybrid IAS, with green corresponding to increments that satisfy the local convexity condition for the inverse gamma, and yellow to the complement.
	The remaining panels of Figure \ref{fig:qr ind} show the corresponding pseudocolor images for the global hybrid algorithm at the switching iteration $\overline{t}$ (center), and at the last iteration of global hybrid IAS (right), respectively.

	\begin{figure}
		\centerline{\includegraphics[width=3.8cm]{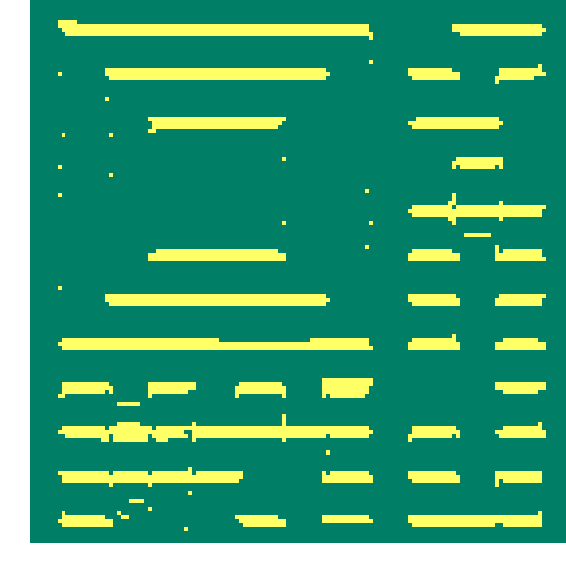} \quad \includegraphics[width=3.8cm]{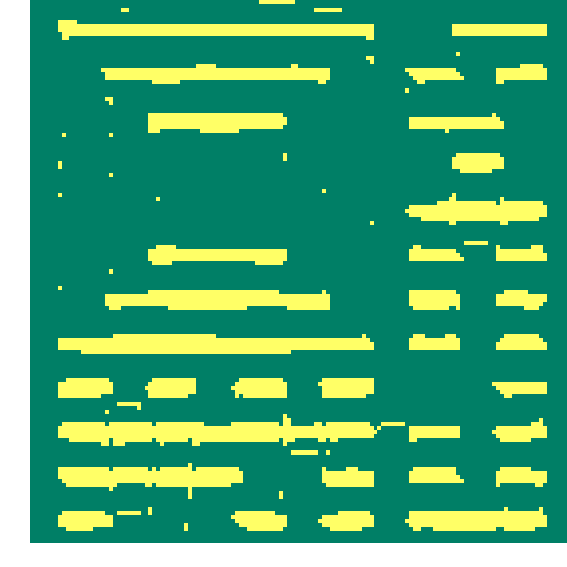} \quad \includegraphics[width=3.8cm]{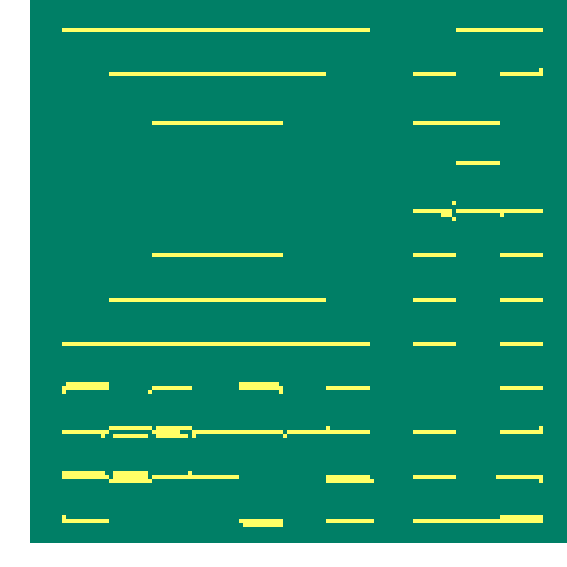}} 
		\centerline{\includegraphics[width=3.8cm]{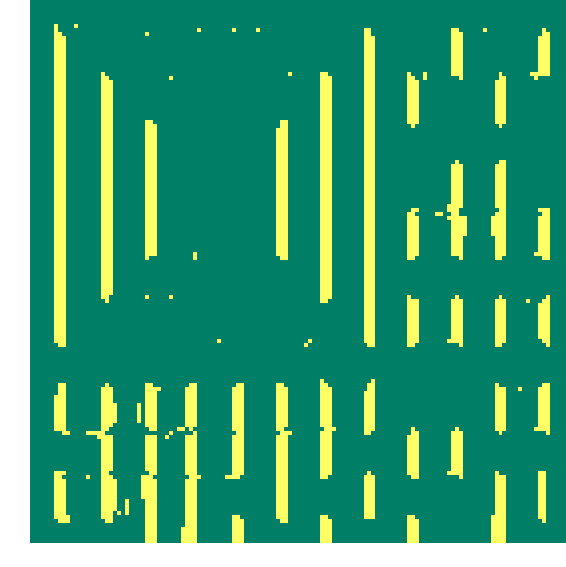} \quad \includegraphics[width=3.8cm]{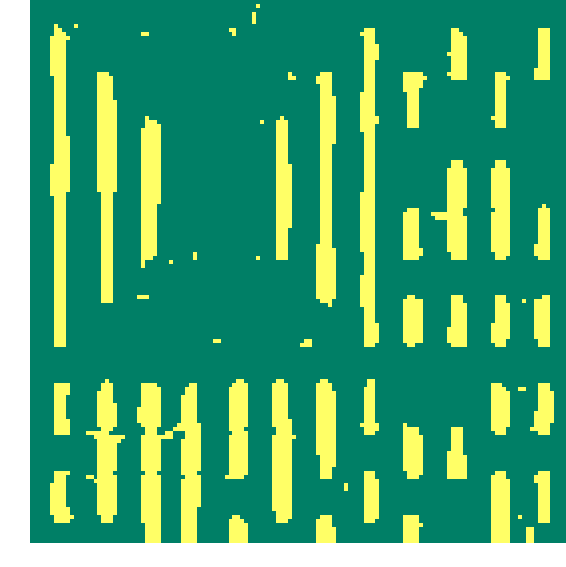} \quad \includegraphics[width=3.8cm]{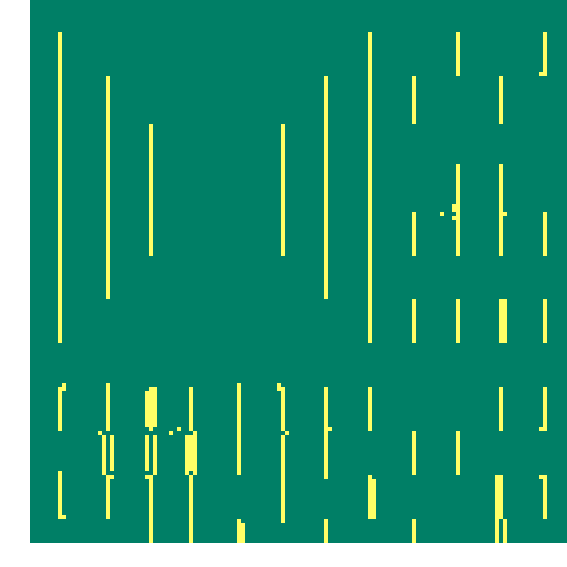}}
		\caption{	\label{fig:qr ind} Pseudocolor image of the indices to the variances $\theta_j$ for vertical (top) and horizontal increments (bottom) with color coding indicating whether $\theta_j< \bar{\theta}$ (green) or $\theta_j \geq \bar{\theta}$ (yellow). 
			The right panels represent the final iteration of the local hybrid algorithm, middle panels the iteration $\overline t-1$ right before the switch of the global hybrid algorithm, and right panels the final iteration of global hybrid algorithm.}
	\end{figure}

	\paragraph{Example 3} In the third example, we consider the problem of estimating a nearly black two-dimensional object. The generative model  is a starry night impulse image, defined as a distribution  on $\Omega =[0,1]\times[0,1]$,
	\[
	d\mu(p) = \sum_{k=1}^J a_k \delta (p -p_k) dp, \quad p_k \sim{\rm Uniform}(\Omega), \quad a_k \sim{\rm Uniform}([1.5,2])\,,
	\]
	is observed through a Gaussian convolution kernel - see (\ref{eq:gauss kern}), with the discrete and noisy data at observation points $q_j\in\Omega$ given by
	\[
	b_j = \int_\Omega A(q_j,p') d\mu(p') + \varepsilon_j = \sum_ {k=1}^K a_k A(q_j,p_k) +\varepsilon_j.
	\] 
	To solve the inverse problem, we subdivide the image  $\Omega$ into $n = 128\times 128  = 16\,384 $ pixels, denoted by $\Omega_\ell$, and let $\mA$ be the matrix representing the discretized kernel,
	
	\[
	\int_\Omega A(q_j,p) d\mu(p) \approx \sum_{\ell =1}^n \underbrace{|\Omega_\ell| A(q_j,q_\ell')}_{=\mA_{j\ell}} x_\ell, \quad x_\ell = \frac 1{|\Omega_\ell|}\int_{\Omega_\ell} d\mu(p),
	\]
	where $q_\ell'$ denotes the center point of the pixel $\Omega_\ell$ and $|\Omega_\ell|$ is its area. The number of observation points is $m = 64\times 64 = 4\,096$ and the noiseless signal is corrupted by scaled white noise with standard deviation approximately 1.8\% of the maximum noiseless signal.  In this case, since the signal itself is sparse, no change of variable is needed. Figure~\ref{fig:starry true} shows the original impulse image characterized by $k=80$ non-zero points, and the noisy blurred image with kernel width $w=0.015$.
	\begin{figure}
		\centerline{\includegraphics[height=4.2cm]{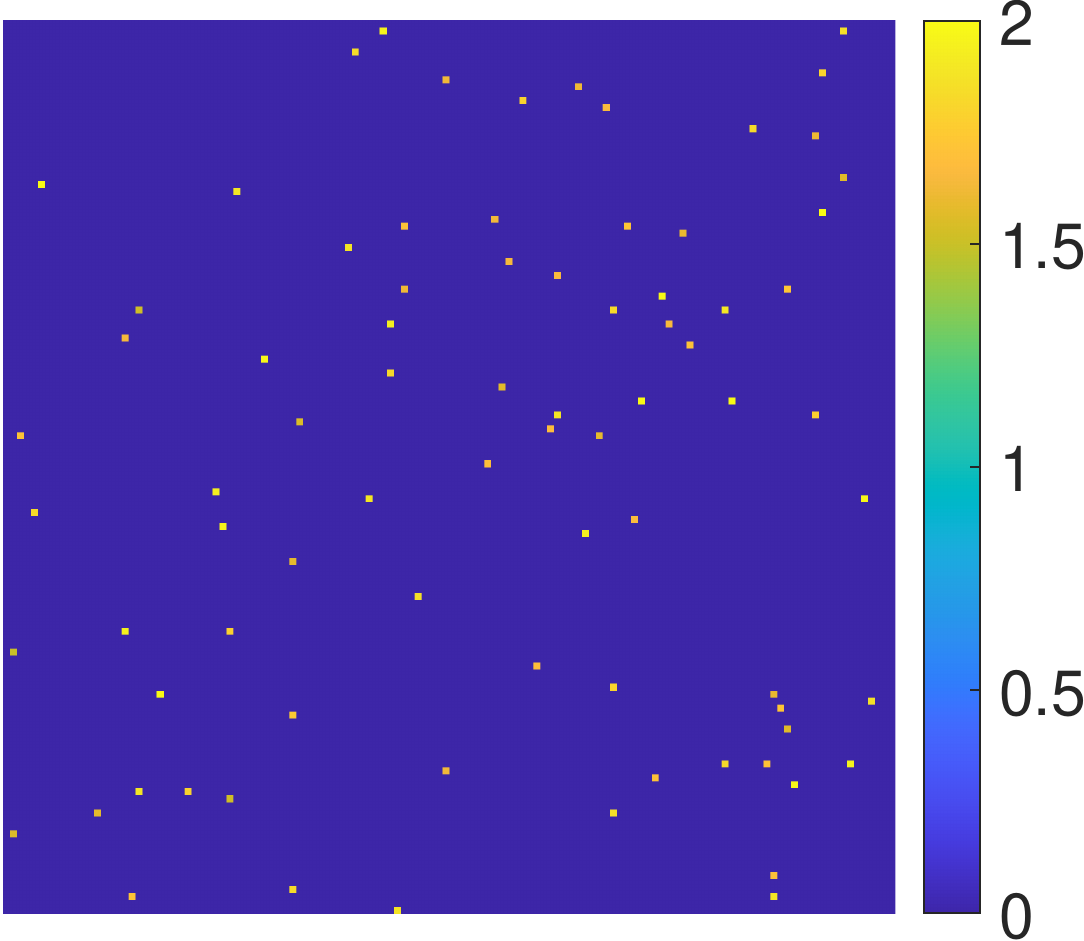} \quad 
			\includegraphics[height=4.1cm]{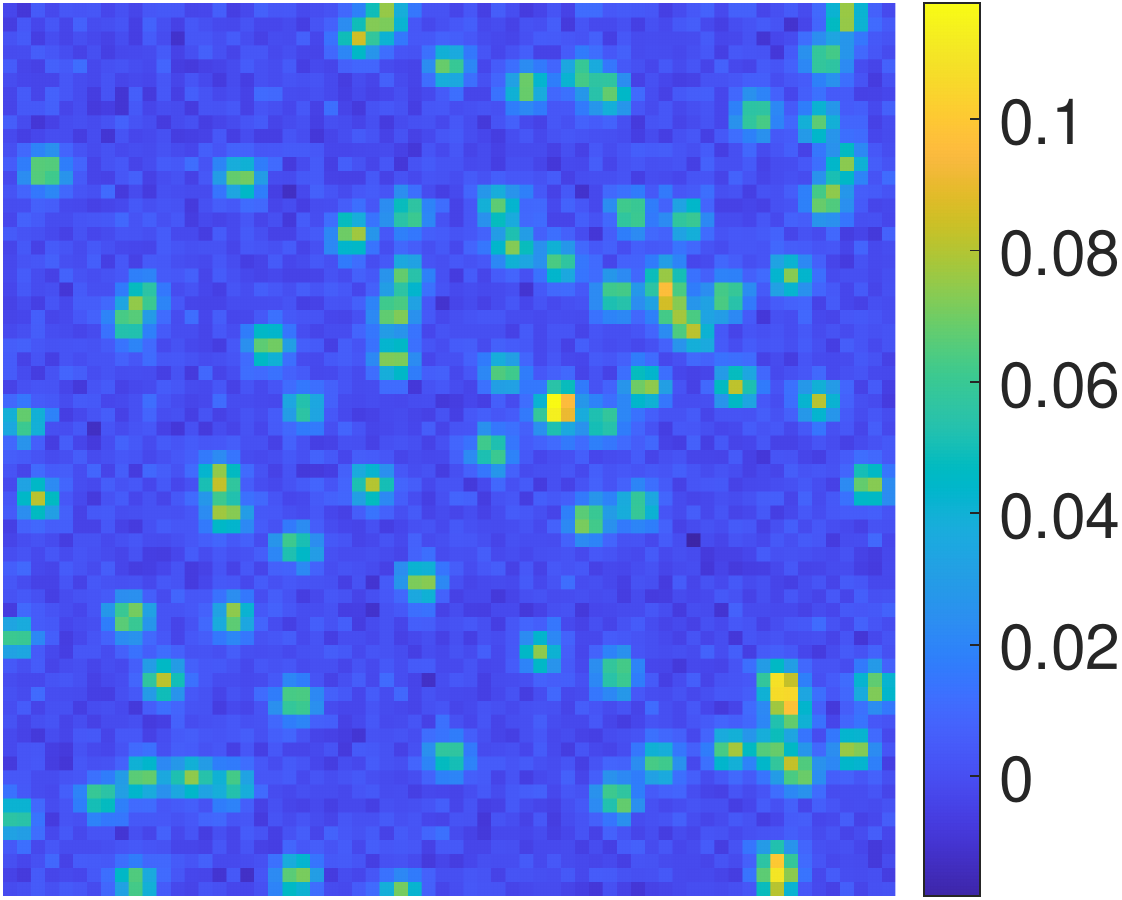}}
		\caption{\label{fig:starry true} Left: Original generative impulse image plotted on a $128 \times 128$ grid as a pixel image. This is the discretization used in the inverse solver, so the pixel image shown here represents the best reconstruction that the algorithm could produce. The reconstructions are compared with this image.  
			Right panel: The $64 \times 64$ blurred and noisy observation, degraded by Gaussian blur and additive while Gaussian noise, scaled so as to achieve SNR $\approx$ 25.}
	\end{figure}

	The restored images and the variances $\theta$ represented as pixel images obtained with the IAS algorithm with gamma and inverse gamma hyperpriors, and the local and global hybrid IAS algorithms   are shown in Figure \ref{fig:starry res}. The differences in the four algorithms are clearly visible from the estimates of the variance $\theta$ and the one-dimensional profiles along the dotted horizontal cut line in the image. The changes in the image of the variances estimated by the IAS with gamma hyperprior (top row, middle) is less pronounced than for the estimates obtained with the other three algorithms, and the intensity of the second star along the cut line (top row, right panel) is significantly lower than in the original image.
	On the other hand, while the reconstruction from the IAS with inverse gamma hyperprior  (second row) is very sharp, the algorithm is too greedy and misses  the second star on the horizontal cut line (right panel). Both hybrid reconstructions (third and fourth row) are sharper than that obtained with the gamma hyperprior image, and reproduce the original profile with higher fidelity than the IAS algorithm with inverse gamma hyperprior.
	
	\begin{figure}
		\centerline{
			\includegraphics[height=3.6cm]{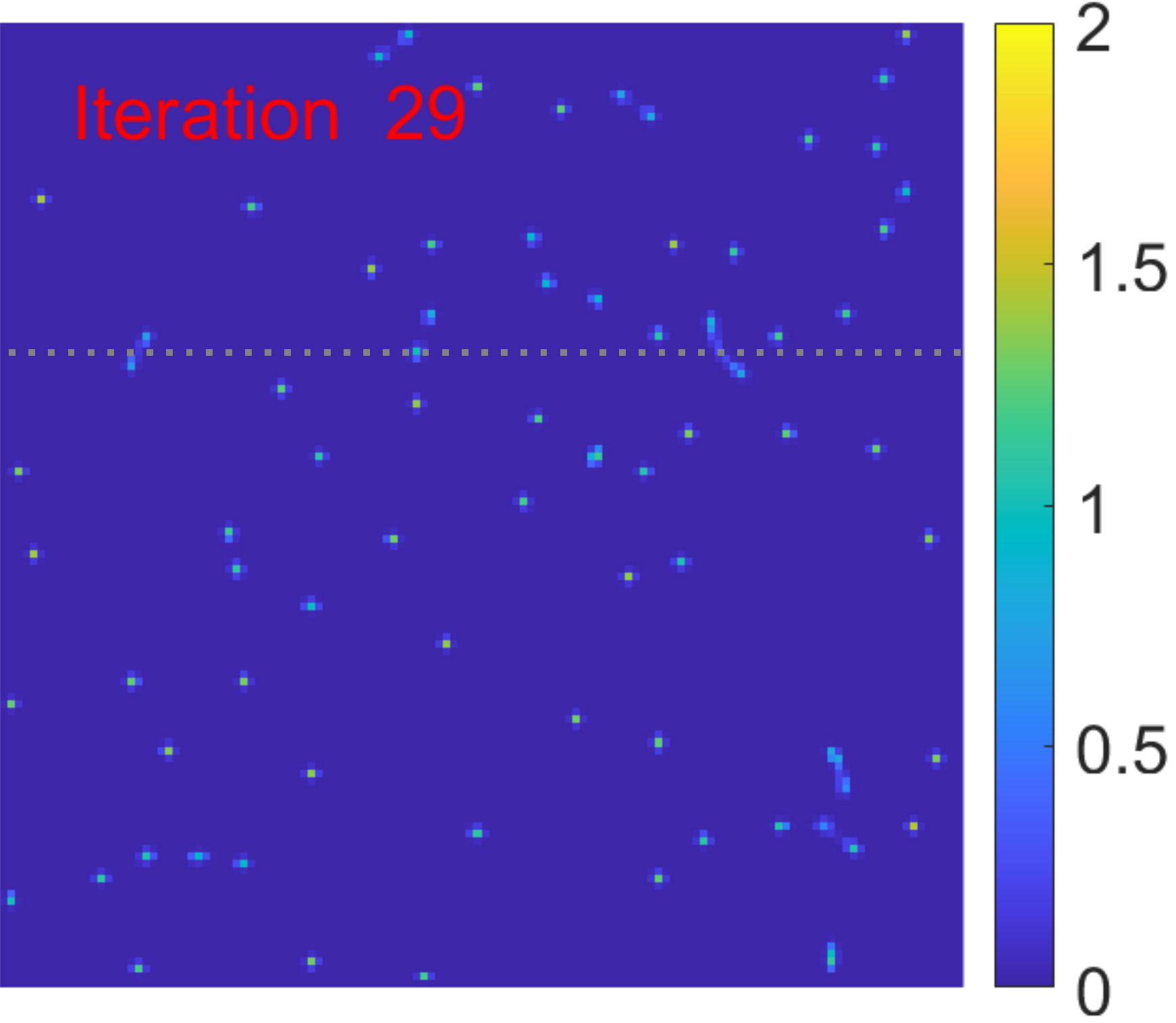} 	
			\quad 
			\includegraphics[height=3.5cm]{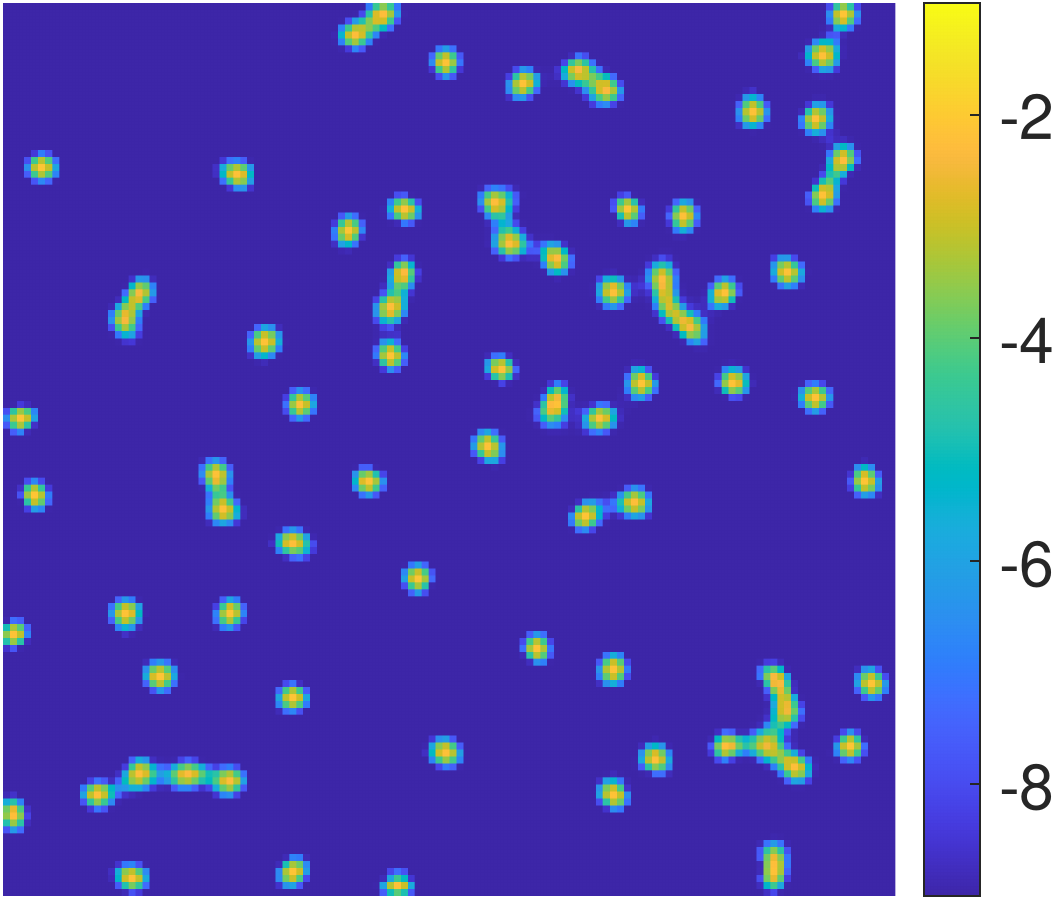}\quad \includegraphics[height=3.5cm]{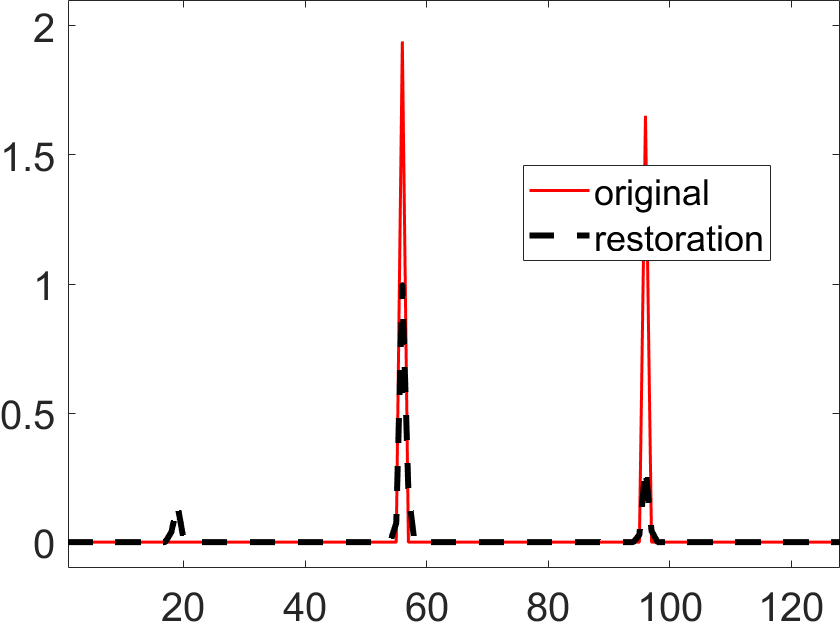}}
		\centerline{
			\includegraphics[height=3.6cm]{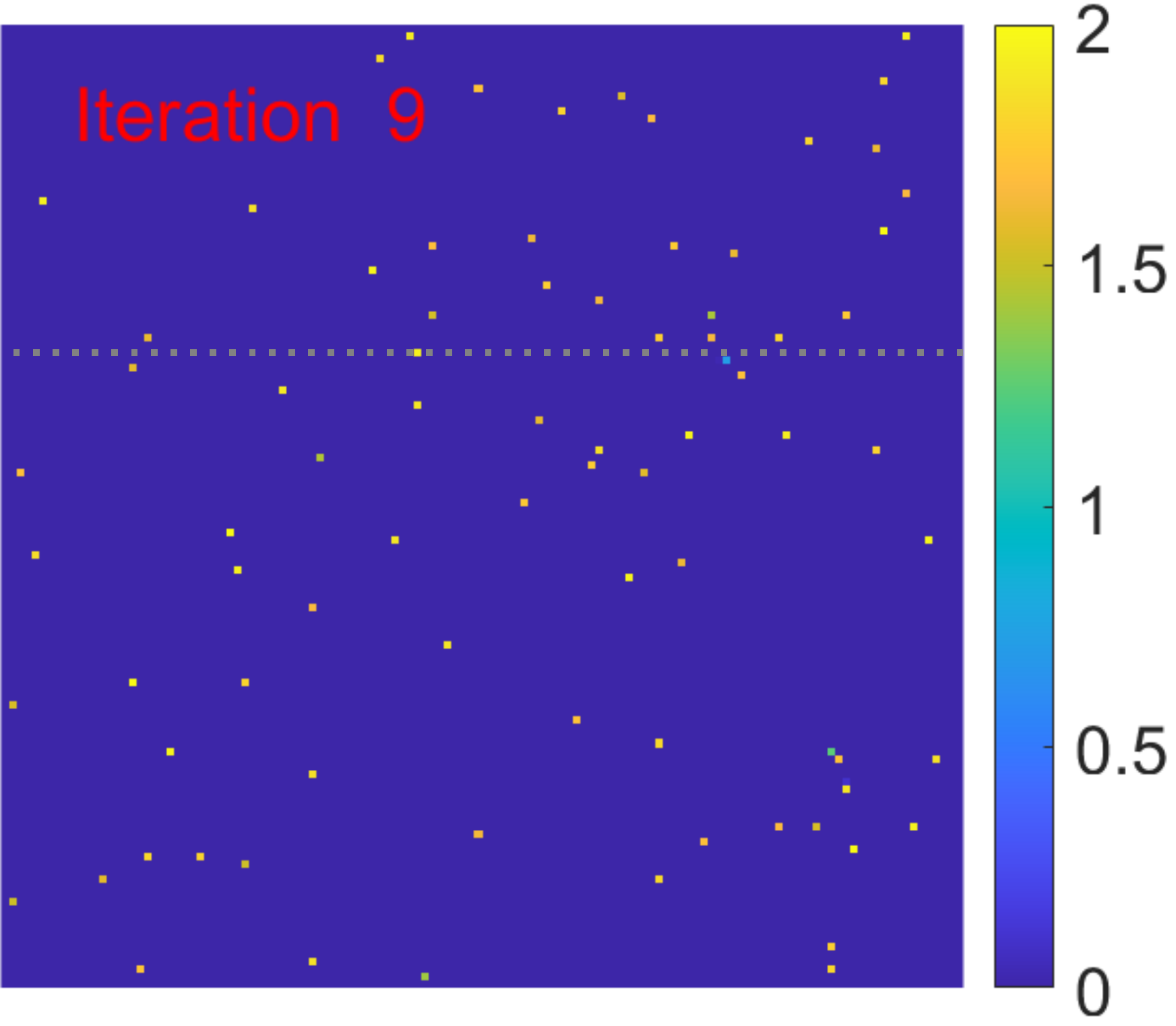} 
			\quad 
			\includegraphics[height=3.5cm]{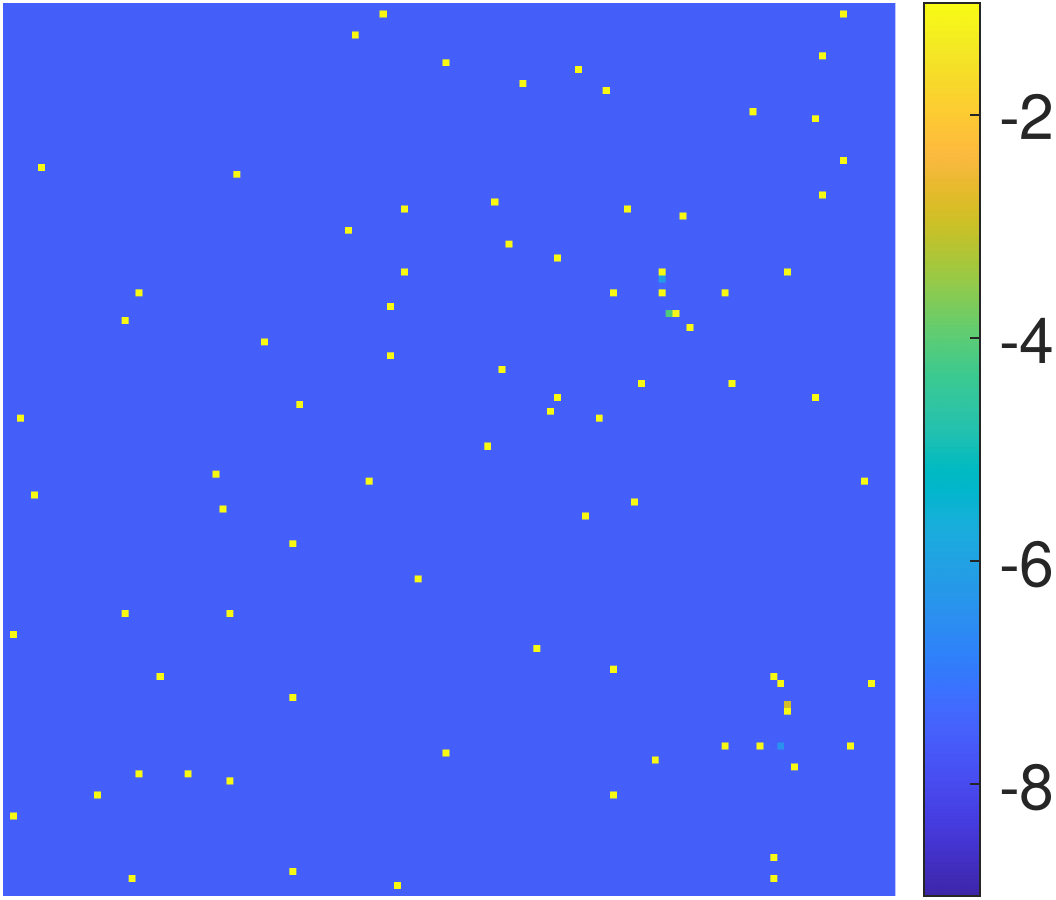}\quad \includegraphics[height=3.5cm]{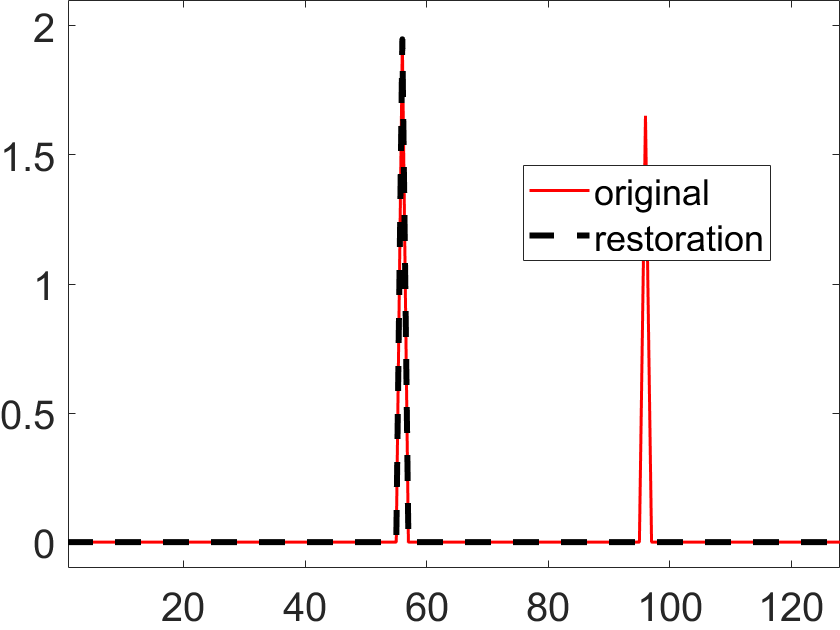}}
		\centerline{
			\includegraphics[height=3.6cm]{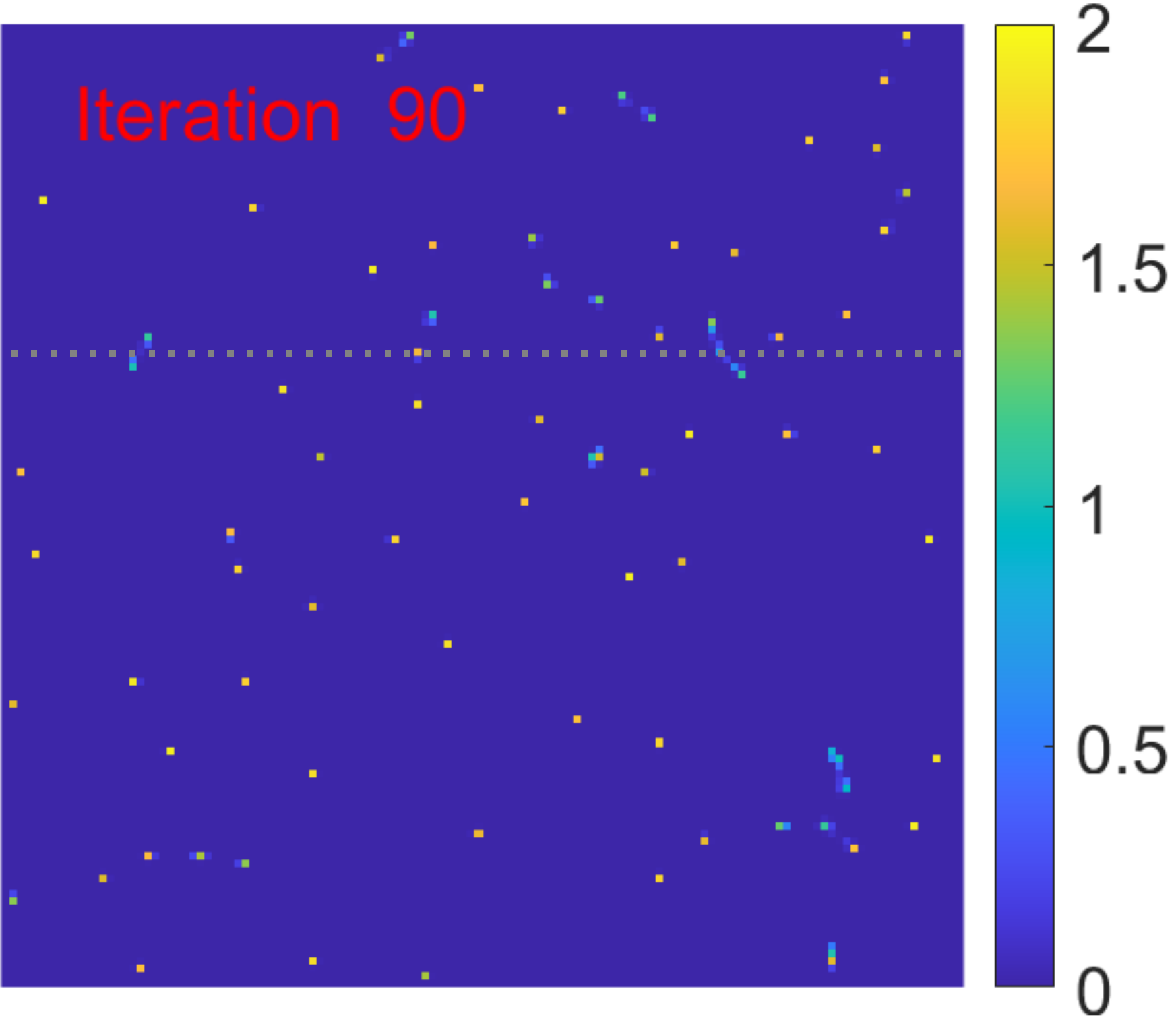}
			\quad 
			\includegraphics[height=3.5cm]{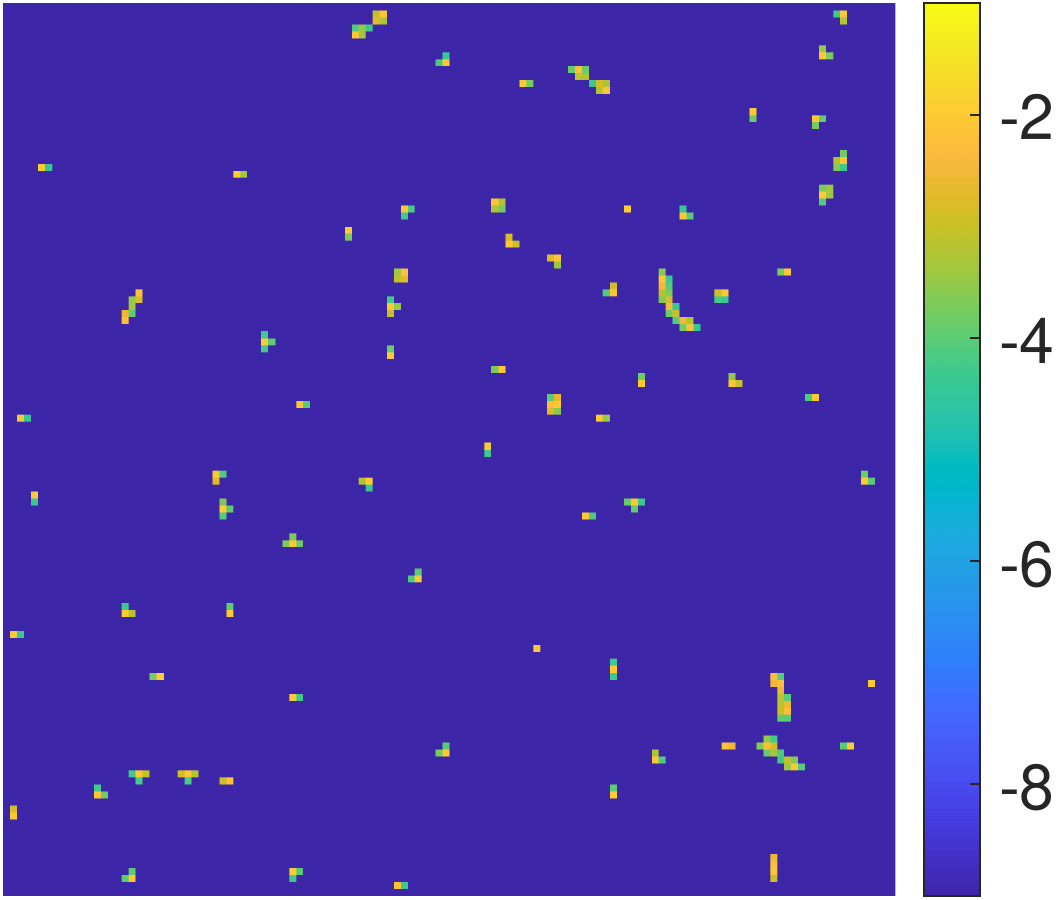}\quad \includegraphics[height=3.5cm]{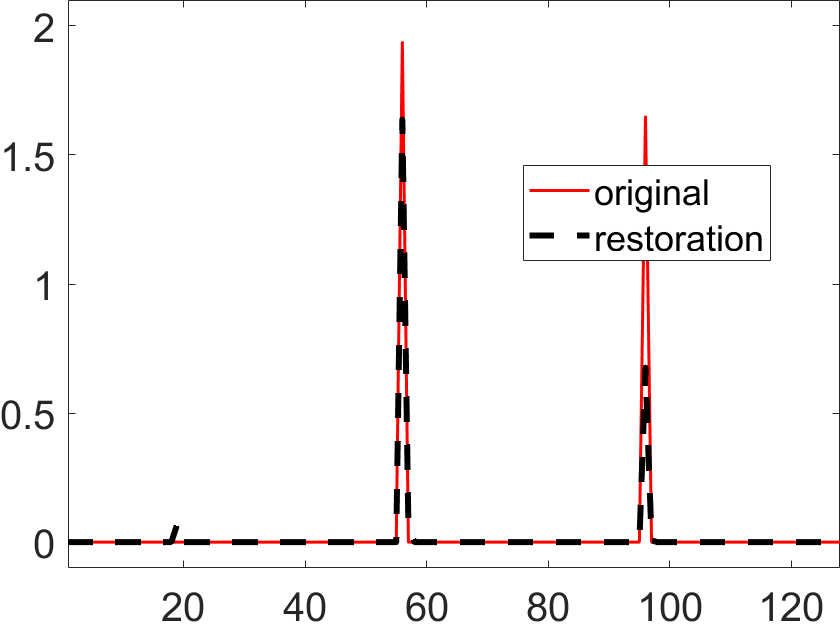}}
		\centerline{
			\includegraphics[height=3.6cm]{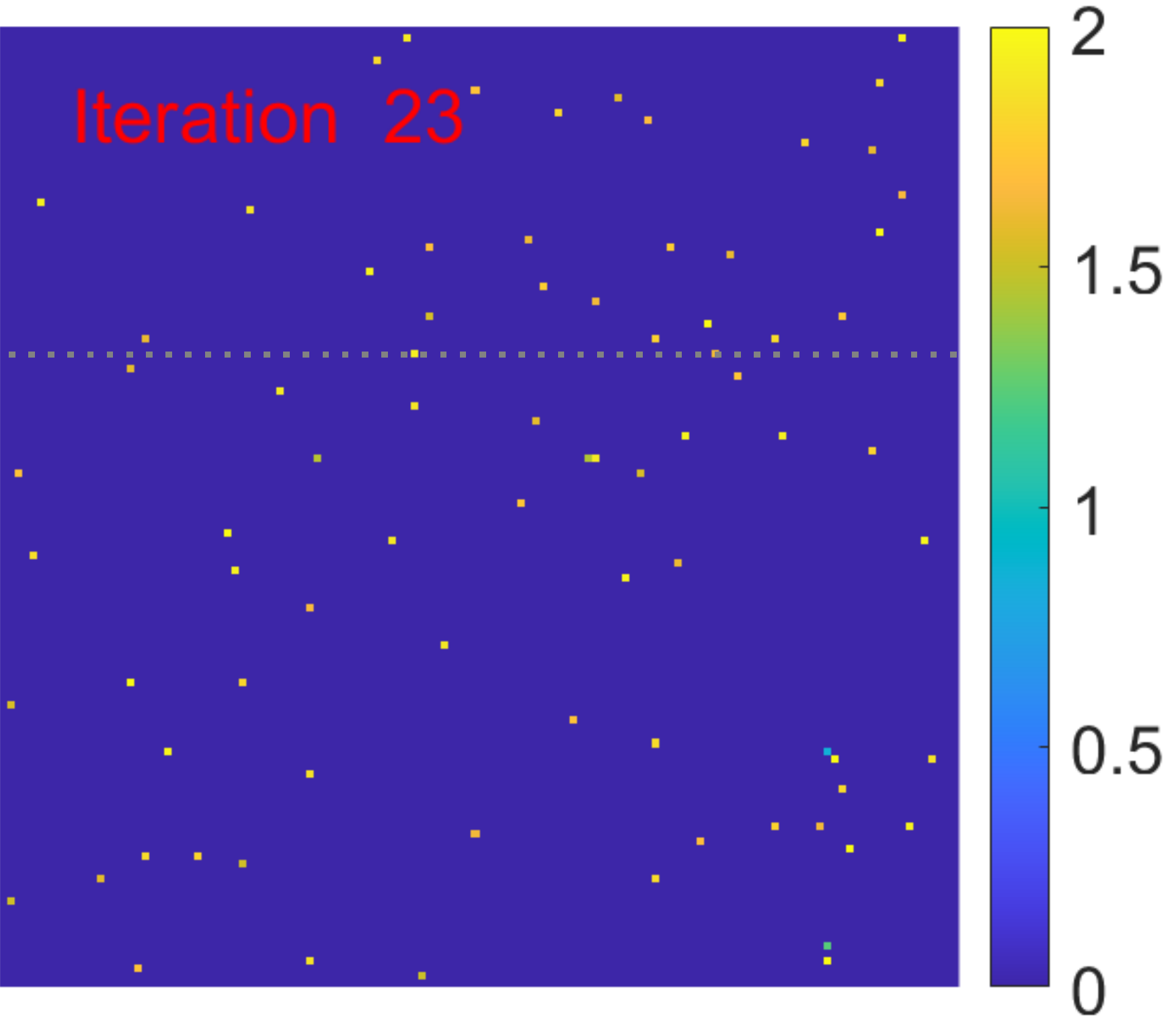}
			\quad 
			\includegraphics[height=3.5cm]{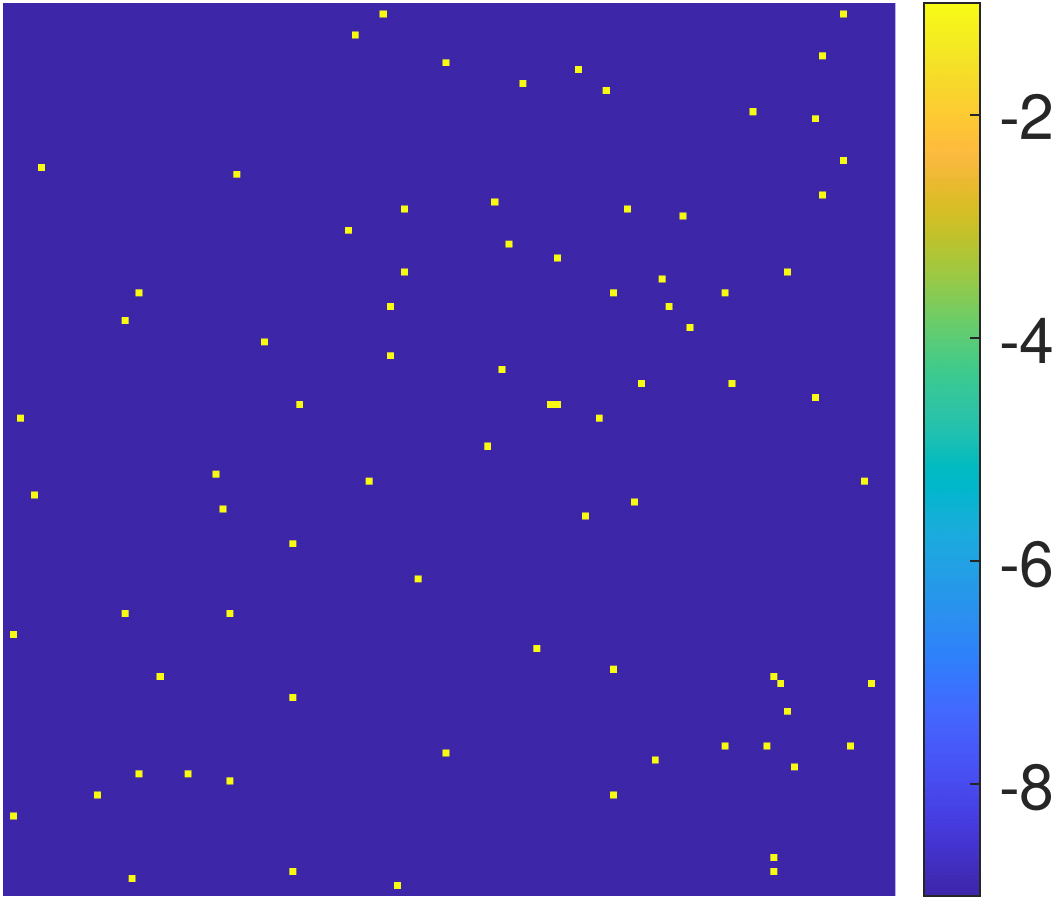}\quad \includegraphics[height=3.5cm]{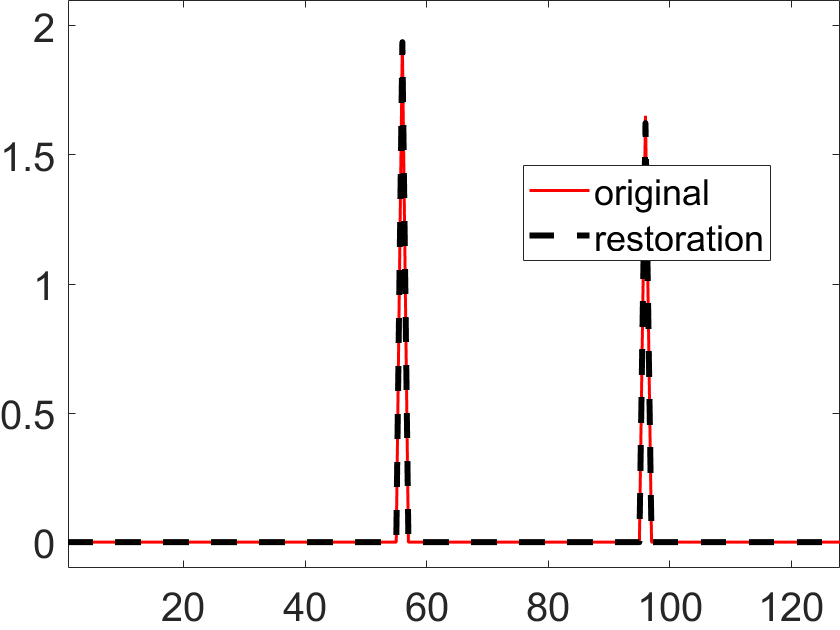}}
		\caption{\label{fig:starry res} Reconstructions, from top to bottom, using gamma and inverse gamma hyperpriors, and the local and global hybrid models (first column). The second column shows the corresponding variances, and the last column shows the reconstructed profiles along the horizontal dotted line across the image. The profiles are compared to the corresponding profile of the generative model represented as a pixel image in the same grid. For the gamma hyperprior in the top row the parameter values are $\eta=10^{-5}$ and $\vartheta=10^{-4}$, for the inverse gamma hypeprior in the second row $\eta=-4.5$ and $\vartheta=10^{-6}$. The hybrid hyperpriors in the bottom rows inherit the parameters from the generative hyperpriors.}
	\end{figure}
	
	Finally, in Figure \ref{fig:starry ind} the behavior of the variances in terms of distribution is shown as a pseudocolor map in the local hybrid case at the final IAS iteration (left panel) and for the global hybrid case  at the switching (middle panel) and final  (right panel) IAS iteration.
	
	\begin{figure}
		\centerline{\includegraphics[height=4.2cm]{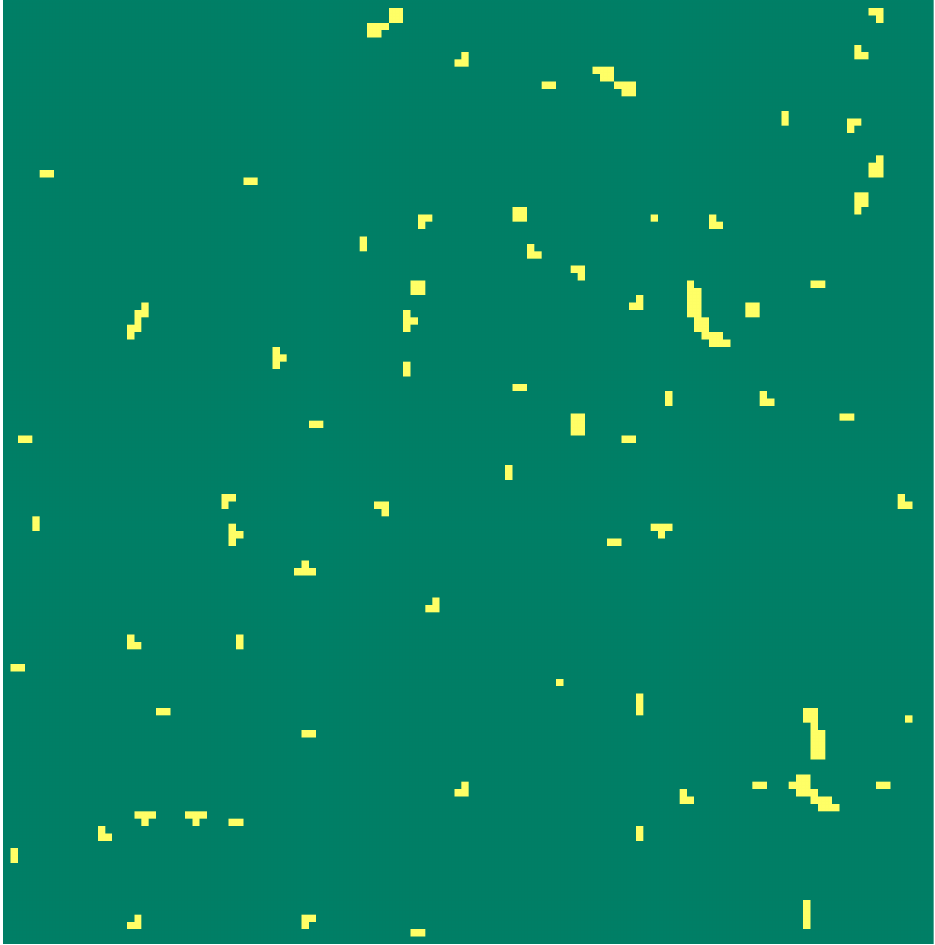}	\quad \includegraphics[height=4.2cm]{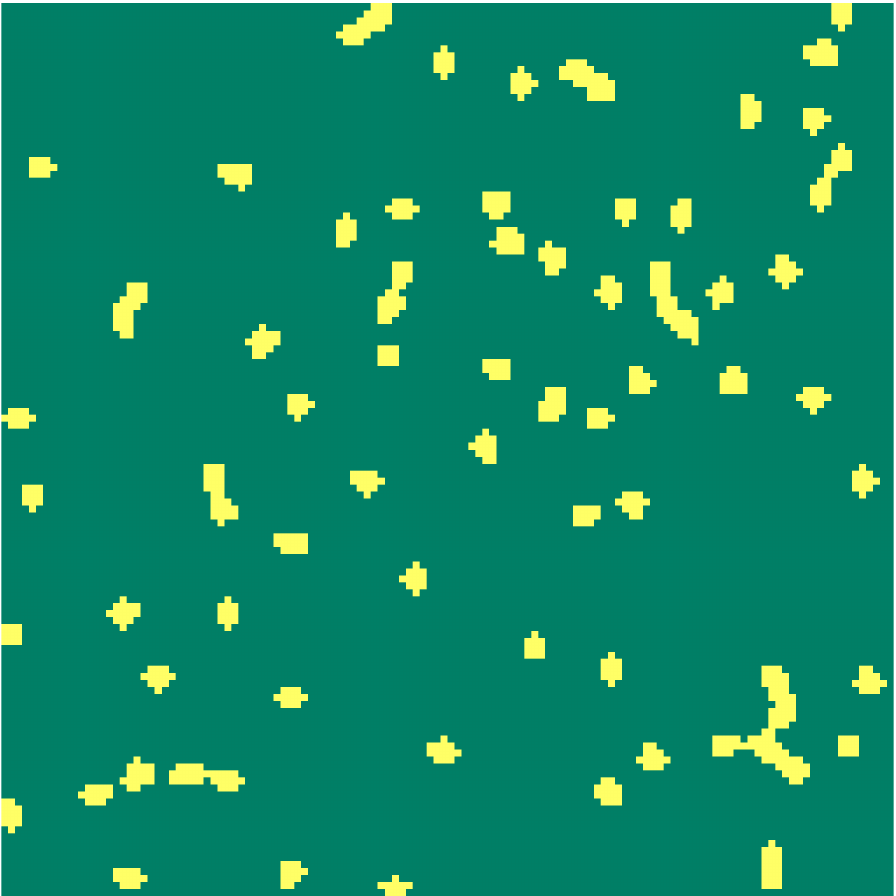} \quad \includegraphics[height=4.2cm]{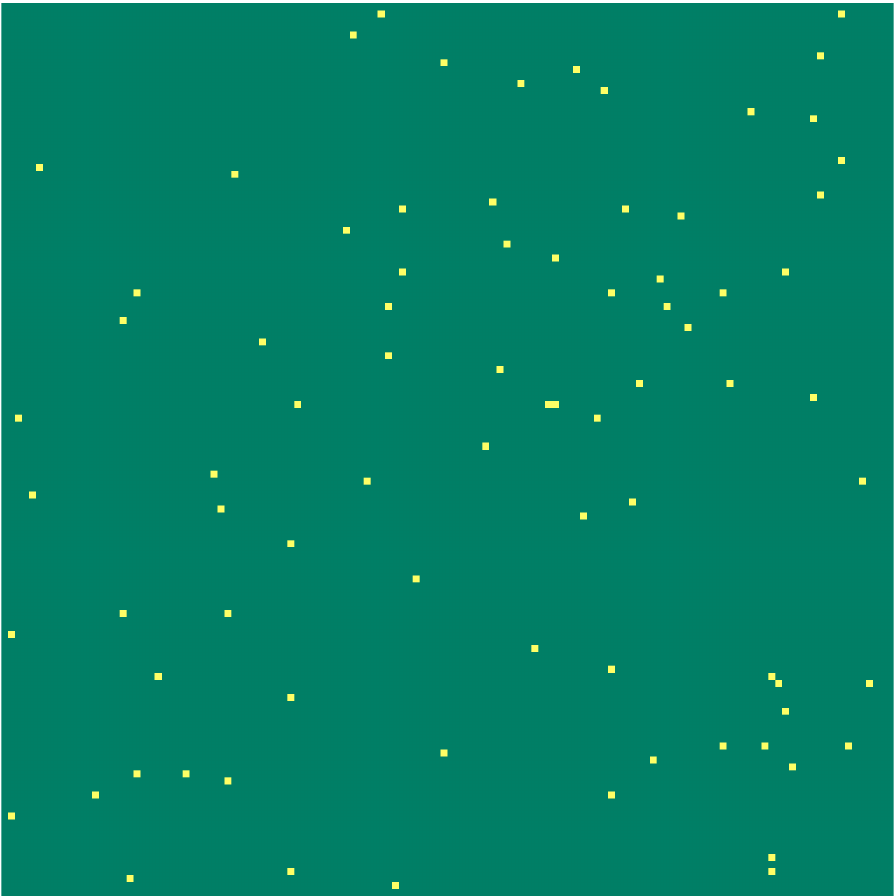}	}
		\caption{\label{fig:starry ind} Image of the variances $\theta_j$ with color coding indicating if $\theta_j< \bar{\theta}$ (green) or $\theta_j \geq \bar{\theta}$ (yellow). Left panel represents the final iteration of the local hybrid algorithm, middle panel the iteration $\bar{t}-1$, right before the switch of the global hybrid algorithm, and right panel the final iteration of global hybrid algorithm.}
	\end{figure}
	
	\section{Conclusions}
	In the present work, we discuss the minimization of conditionally Gaussian hypermodels under the adoption of generalized gamma hyperpriors. Based on the results derived in \cite{CPrSS}, the two proposed hybrid algorithms, namely the local and global hybrid IAS, exploit the global convexity ensured by gamma hyperpriors ($r=1$) and the stronger sparsity promotion of the generalized gamma hyperpriors with $r<1$. The local hybrid hypermodel preserves the global convexity characterizing the gamma hyperpriors and, as confirmed by numerical examples, is particularly effective in cleaning the background, while not ensuring a sharp recovery of sudden discontinuities in the signals. On the other hand, the global hybrid hypermodel, which relies on the detection of a suitable initial guess for the minimization of the locally convex hypermodel ${\mathscr M}_2$, returns high quality restorations at the expense of global convexity.

	\section*{Acknowlwdgements}
	
	The work of DC was partly supported by the NSF grant DMS-1522334, and of ES by the NSF grant DMS-1714617. Part of this work was done while the authors DC, MP and ES were visiting Institut Henri Poincar\'{e} during the worshop ``The Mathematics of Imaging'' in March - April 2019 and the Institue for Mathematics and its Applications during the workshop ``Computational Imaging'' in October 2019. The hospitality and support of IHP and IMA are gratefully acknowledged.

\end{document}